\documentclass[12pt]{article}
\usepackage{amsthm,amsmath,amsfonts,amssymb}
\usepackage[round]{natbib}
\usepackage{color}
\usepackage[colorlinks,allcolors=blue]{hyperref}
\usepackage{url}
\usepackage{graphicx}
\usepackage[T1]{fontenc}
\usepackage[titletoc,title]{appendix}
\usepackage[left=2.8cm,right=2.8cm,top=2.5cm,bottom=2.5cm]{geometry}
\nocite{*}
\newtheorem{definition}{Definition}
\newtheorem{theorem}{Theorem}
\newtheorem{lemma}{Lemma}
\newtheorem{remark}{Remark}
\newtheorem{corollary}{Corollary}
\newtheorem{example}{Example}

\title{\sc{\large{\textbf{Bayesian and minimax estimators of loss}}}}
\author{{\sc{\normalsize
Christine Allard$^{a}$ \ \& \
{\'E}ric Marchand$^{a}$}}\\
\centerline{\small{$^{a}$Universit\'e de Sherbrooke, D\'epartement de math\'ematiques, Sherbrooke, Qc, CANADA, J1K 2R1}}\\
\centerline{\small{email: \texttt{christine.allard@usherbrooke.ca}; \texttt{eric.marchand@usherbrooke.ca}}}\\
}
\date{}
\providecommand{\keywords}[1]{\textbf{\textit{Keywords---}} #1}
\begin{document}
\maketitle
\begin{abstract}
We study the problem of loss estimation that involves for an observable $X \sim f_{\theta}$ the choice of a first-stage estimator $\hat{\gamma}$ of $\gamma(\theta)$, incurred loss $L=L(\theta, \hat{\gamma})$, and the choice of a second-stage estimator $\hat{L}$ of $L$.   We consider both:  (i) a sequential
version where the first-stage estimate and loss are fixed and optimization is performed at the second-stage level, and (ii) a simultaneous version with a Rukhin-type loss function designed for the evaluation of $(\hat{\gamma}, \hat{L})$ as an estimator of $(\gamma, L)$.

We explore various Bayesian solutions and provide minimax estimators for both situations (i) and (ii).   The analysis is carried out for several probability models, including multivariate normal models $N_d(\theta, \sigma^2 I_d)$  with both known and unknown $\sigma^2$, Gamma, univariate and multivariate Poisson, and negative binomial models, and relates to different choices of the first-stage and second-stage losses.  The minimax findings are achieved by identifying least favourable of sequence of priors and depend critically on particular Bayesian solution properties, namely situations where the second-stage estimator $\hat{L}(x)$ is constant as a function of $x$.
\end{abstract}
\section{Introduction}
Reporting on the precision of statistical decisions, whether it relates to a standard error of an estimate in multiple regression or survey sampling, the power of a test, or the coverage probability of an interval estimate etc., are central to the practice of statistics.
Whereas a frequentist approach typically prescribes the level of risk prior to the collection of data, while a Bayesian approach typically relates to post-data inference to assess precision, such correspondences are not exclusively the case and various approaches have been presented in the literature, for instance by \cite{Berger1985}; \cite{Goutis1995}, as well as the references therein.
The setting or search for efficient accuracy reports thus matters and we investigate here various issues and optimality properties cast in a loss estimation framework.
The seminal work of \cite{Johnstone1988} put forth a framework for loss estimation with an emphasis on the multivariate normal model and the discovery of improvements in terms of frequentist risk on usual unbiased estimators.  Earlier work in this direction includes that of \cite{Sandved1968}, but many researchers engaged after Johnstone's work with further investigation towards various extensions in terms of different contexts and other models (e.g., multivariable regression, model selection, spherical symmetry), optimality properties (e.g., admissibility and dominance), and further related issues (e.g.,   \cite{Boisbunon2014};  \cite{Berger1995}; \cite{FourdrinierStraw2003}; \cite{FW1995a, FW1995b, FW2012};  \cite{FourdrinierL2008}; \cite{Lele1992, Lele1993}; \cite{LuBerger1989}; \cite{Maruyama1997}; \cite{Matsuda2019}; \cite{Matsuda2022}; \cite{Narayanan2015}; \cite{WanZou2014}).

The approach which we adopt and study is that of accompanying an estimator $\hat{\gamma}$ of an unknown parameter $\gamma(\theta)$ by an estimator $\hat{L}$ of a first-stage loss $L=L(\theta, \hat{\gamma}(x))$.  To assess the accuracy of such a $\hat{L}$, we work with a second-stage loss typically of the form
$W(L,\hat{L})$.   We explore Bayesian estimators $\hat{L}_{\pi}$ of $L$ and their properties for a given prior density $\pi$.   We also address minimax optimality for estimating $L$ and provide minimax solutions, which capitalize on the behaviour of Bayesian solutions.  Minimax solutions serve as a benchmark for evaluating competing estimators, and  have been relatively unexplored in the context of loss estimation.

A particular interesting and different approach, which combines both the first-stage estimation process and the second-stage loss estimation component  was proposed and analyzed by \cite{Rukhin1987, Rukhin1988}.    Here again, we provide minimax solutions $(\hat{\gamma}, \hat{L})$ which depend critically on the existence of prior densities and associated Bayesian estimators $\hat{L}_{\pi}$ that do not depend on the observed data.

The paper is organized as follows.    Section 2 describes the adopted language of loss estimation namely with the aid of definitions and notations.  Section 3 explores Bayesian solutions in cases where both the first-stage and second-stage estimators $\hat{\gamma}_{\pi}$ and $\hat{L}_{\pi}$ are derived with respect to the same prior, with an emphasis on varied choices of the first-stage and second-stage losses, and namely departures from the ubiquitous squared error second-stage loss in the literature.   Examples of models include normal with or without a known covariance structure, Gamma, univariate and multivariate Poisson, Negative Binomial, and location exponential.  We also come across a surprising large number of cases where the Bayes estimator $\hat{L}$ of loss is constant as a function of the observable data $x$, including cases where even the posterior distribution $L|x$ is free of $x$, and explore links between Bayesian  and unbiased estimators of loss.

Section 4 provides minimax estimators $\hat{L}$ of loss $L(\theta, \hat{\gamma}(x))$ for different combinations of probability models and losses, and when the first-stage estimate has previously been obtained.  To the best of our knowledge, the only known previously analyzed case involves normal models and first and second-stage squared error losses \cite{Johnstone1988}.  We extend the finding to a wider class of first and second-stage losses, to an unknown covariance structure, and proceed with minimax results for Gamma models.   The results are obtained through the determination of a least favourable sequence $\{\pi_n\}$ of prior densities and an extended Bayes estimator with constant risk.

 Finally, we obtain minimax solutions $(\hat{\gamma}_0, \hat{L}_0)$ for estimating
$(\gamma, L)$ simultaneously under Rukhin-type losses, by using a sequence of priors approach to obtain an extended Bayes estimator with constant frequentist risk.  Properties established or observed in Section 3, namely the constancy of Bayesian solutions $\hat{L}_{\pi}(x)$ as functions of $x$ play a key role for the analysis.
\section{Definitons and other preliminaries}

Throughout, we consider a model density $f_{\theta}$, $\theta \in \Theta$, for an observable $X$, a parametric function $\gamma(\theta)$ of interest; mostly taken as identity;  and the loss $L=L(\theta, \hat{\gamma})$ measuring the level of accuracy (or inaccuracy) of $\hat{\gamma}(X)$ as an estimator of $\gamma(\theta)$.  Such a loss $L$ is referred to as a first-stage loss, while loss $W(L,\hat{L}) \in [0,\infty)$ used for estimating $L$ by $\hat{L}(X)$ is referred to as the second-stage loss.

\begin{remark}
\label{remarklosschoice}   A common and default choice in the literature for the second-stage loss $W$ has been squared error $(\hat{L}-L)^2$ loss.   This has been described as a matter of convenience, but it is also the case that the developments for multivariate normal models, as well as spherically symmetric models, bring into play Stein's lemma
and convenient loss estimation representations arising for squared error loss (e.g. \cite{FSW2018}.  However, since a loss $L$ is more synonymous with a positive scale parameter than a location parameter, it seems desirable to consider a scale invariant second-stage loss of the form $\rho(\frac{\hat{L}}{L})$, plausible choices satisfying the bowl-shaped property $\rho(t)$ decreasing for $t \in (0,1)$ and increasing for $t>1$.   The developments in this manuscript thus address such losses amongst a wider choice of second-stage losses. Such examples include weighted squared error loss with $\rho(t)=(t-1)^2$ , squared log error loss with $\rho(t)\, = \, (\log t)^2$, symmetric versions with $\rho(t)=\rho(\frac{1}{t})$ like $\rho(t)=t+ \frac{1}{t}-2$ (e.g., \cite{Mozgunov2019}), entropy loss with $\rho(t)\,=\, t -  \log t - 1$, and variants of the above (except log error) with $\rho_m(t) = \rho(t^m)$, $m  \neq 0$ and the case $m = -1$ being most prominent for squared error and entropy.
\end{remark}

For a given prior density $\pi$ for $\theta$ and estimate $\hat{\gamma}(x)$, a
Bayes estimator $\hat{L}_{\pi}(X)$ of $L=L\big(\theta, \hat{\gamma}(x)\big)$ minimizes in $\hat{L}$ the expected posterior loss $\mathbb{E}\big(W(L,\hat{L})|x   \big)$.  As addressed in Section \ref{sectionbayesian}, it is particularly interesting  to study cases where $\hat{\gamma} \equiv \hat{\gamma}_{\pi}$, i.e., the Bayes estimator of $\gamma(\theta)$ under the same prior.   However, it is still useful to consider the more general context, for instance because the first-stage estimator may be imposed and not be Bayesian, or a theoretical assessment of a least favourable sequence of priors such as the one pursued in Section \ref{subsectionminimax} requires it.  In such a framework, the stated goal is how to report on sensible estimates of $L$ for a given $\hat{\gamma}(x)$.

The frequentist risk performance of a given estimator $\hat{L}$ for estimating a loss $L\,=\, L(\theta, \hat{\gamma}(x))$ is given by
\begin{equation}
\label{risk}
R(\theta, \hat{L}) \, = \, \mathbb{E}_{\theta} \big( W(L, \hat{L}(X))  \big)\,, \theta \in \Theta\,,
\end{equation}
and different choices of $\hat{L}$ can be compared leading to the usual definitions of dominance, admissibility, and inadmissibility.  Our findings  relate to the minimax criterion with estimator $\hat{L}_m(X)$ being a minimax estimator of $L$ whenever  $\sup_{\theta \in \Theta} \{R(\theta, \hat{L}_m)\} \, = \,  \inf_{\hat{L}} \sup_{\theta \in \Theta} \{R(\theta, \hat{L})\}$.

Another criterion present in the literature, that sometimes interacts with Bayesianity, is that of unbiasedness.
For a given estimator $\hat{\gamma}(X)$ of $\gamma(\theta)$, an estimator $\hat{L}(X)$ of loss $L(\theta, \hat{\gamma})$ is said to be unbiased if
\begin{equation}
\label{definitionunbiased}
\mathbb{E}_{\theta} \, \hat{L}(X) \, = \, \mathbb{E}_{\theta} L\big(\theta, \hat{\gamma}(X)\big), \hbox{ for all } \theta \,,
\end{equation}
i.e.,  $\hat{L}(X)$ is an unbiased estimator of the frequentist risk of $\hat{\gamma}(X)$ as an estimator of $\gamma(\theta)$.
\begin{example}
\label{exampleJohnstone}
As an illustration, consider the normal model $X \sim N_d(\theta, \sigma^2 I_d)$ with known $\sigma^2$, the benchmark estimator $\theta_0(X)=X$ of $\theta$, and the incurred squared error loss  $L=L(\theta, \hat{\theta}(X))\,=\,\|X-\theta\|^2$.
Johnstone (1988) showed that the estimator $\hat{L}_0(X) \, = \, d \sigma^2$; equal here to the mean squared error risk of $\hat{\theta}_0$ making it an unbiased estimator of loss;  matches the (generalized) Bayesian estimator $\hat{L}_{\pi_0}(X)$ of $L$ for second-stage squared error loss $W(L,\hat{L})\,=\, (\hat{L}-L)^2$ and the uniform prior density $\pi_0(\theta)=1$.    He also established that $\hat{L}_{\pi_0}$ is minimax for all $d \geq 1$, admissible for $d \leq 4$, and inadmissible for $d \geq 5$ providing dominating estimators $\hat{L}$ in such a latter case.   Such dominating procedures are necessarily minimax.   Our findings (Section \ref{subsectionnormalmodels}) for normal models relate to minimaxity for different choices of first-stage ($L$) and second-stage ($W$) losses, and address the case of unknown $\sigma^2$ (Section \ref{sectionnormalunknownsigma^2}).
\end{example}

Rukhin (1987, 1988) studied the efficiency of estimators, namely in terms of admissibility,  with his proposal to combine the two stages of estimation and to measure the performance of the pair $(\hat{\gamma}, \hat{L})$ for estimating $\big(\gamma(\theta), L  \big)$, with $L=L(\theta, \hat{\gamma})$ by the loss
\begin{equation}
\label{lossrukhin}
\hat{L}^{-1/2} L(\theta, \hat{\gamma}) + \hat{L}^{1/2}\,,
\end{equation}
as well as extensions
\begin{equation}
\label{lossrukhingeneral}
W(\theta, \hat{\gamma}, \hat{L}) \, = \,  h'(\hat{L}) \, L(\theta, \hat{\gamma}) - h'(\hat{L}) \, \hat{L} \, + \,h(\hat{L})\,,
\end{equation}
$h$ being an increasing and concave function on $(0,\infty)$, and the former being a particular case of the latter for $h(\hat{L})\,=\, 2 \hat{L}^{1/2}$.  The two components are referred to as error of estimation and precision of estimation, and such a loss is appealing namely since $\hat{L}=L$ minimizes the loss for fixed $L$ and since the Bayesian estimator of $L$ for a given prior $\pi$ is equal to the posterior expectation $\hat{L}_{\pi}(x) \, = \, \mathbb{E}(L|x)$.

\section{Bayesian estimators}
\label{sectionbayesian}
In this section, we record various interesting scenarios concerning Bayesian inference about a given loss $L=L(\theta, \hat{\gamma})$ incurred by estimator $\hat{\gamma}(X)$ for estimating $\gamma(\theta)$.  Different choices of $L$ and the second-stage loss $W(L, \hat{L})$ are considered for the determination of a point estimator $\hat{L}_{\pi}$, and we also describe directly the posterior distribution $L|x$ in some instances.

Section \ref{sectionpairs} deals with situations where a same prior $\pi$ is used to determine both the choices of $\hat{\gamma}_{\pi}$ and $\hat{L}_{\pi}$, with a particular focus on Poisson and negative binomial models.  While these last two examples involve cases where point estimates $\hat{L}_{\pi}(x)$ do not depend on $x$,  Section \ref{secionposteriorfreeofx} deals with specific situations where the posterior distribution $L|x$ does not depend on $x$.    The features and representations provided in the Section will prove useful for the minimax findings of Sections \ref{sectionminimax} and \ref{sectionrukhin}.

\subsection{Examples of $(\hat{\gamma}_{\pi}, \hat{L}_{\pi})$ pairs}
\label{sectionpairs}

From a Bayesian perspective with the same prior $\pi$ as for the first-stage, one could naturally consider the posterior distribution of $L|x$ for inference about $L$.  Minimizing the expected posterior loss $W(L,\hat{L})$ in $\hat{L}$ produces Bayes estimate $\hat{L}_{\pi}$ and the ensemble produces pairs $\big(\hat{\gamma}_{\pi}, \hat{L}_{\pi}\big)$ which we investigate and illustrate in this section.

We begin with the familiar case of squared error loss $L(\theta, \hat{\theta}) = \|\delta-\theta\|^2$ for estimating $\gamma(\theta)=\theta \in \mathbb{R}^d$ based on $X \sim f_{\theta}(\cdot)$.  Assuming that the posterior covariance matrix $Cov(\theta|x)$ of $\theta$ exists, we have $$\hat{\gamma}_{\pi}(X)=\mathbb{E}(\theta|x) \hbox{ with incurred loss } L \, = \, \|\theta - \mathbb{E}(\theta|x)\|^2\,.$$   Now, if the second-stage loss is again squared error loss, i.e., $W(L,\hat{L})=(\hat{L}-L)^2$, then we obtain
\begin{equation}
\label{lhatpi}
  \hat{L}_{\pi}(x)=\mathbb{E}(L|x) \, = \, tr \, Cov(\theta|x)\,.
\end{equation}

\begin{example}
\label{examplelocationfamily}
For location family densities $f_{\theta}(x)=f_0(x-\theta)$, $x, \theta  \in  \mathbb{R}^d$, such that $\mathbb{E}_{\theta}(X)=\theta$, which include spherically symmetric densities $g(\|x-\theta\|^2)$, and non-informative prior $\pi(\theta)=1$, we obtain $\hat{\gamma}_{\pi}(x)\,=\,x $.  Since $x-\theta|x =^d X-\theta|\theta$ for all $x, \theta$, it follows that the posterior distribution $L|x$, which of that of the radius $\|\theta -x\|\big|\theta$, is independent of $x$.   Consequently,  $\hat{L}_{\pi}(x)$ is a constant in terms of $x$.  This observation includes the familiar multivariate normal case with $X|\theta \sim N_d(\theta, \sigma^2 I_d)$ where $L|x \sim \sigma^2 \chi^2_d$ and  $\hat{L}_{\pi}(x) \, = \, d \sigma^2$.  With the posterior distribution of $L$ independent of $x$, Bayes estimators $\hat{L}_{\pi}(X)$ associated with other second-stage losses or even credibility intervals for $L$ will necessarily be independent of $x$.  These last features fit into a more general structure expanded upon with Theorem \ref{theoremlocation} and include proper priors as with the following example.
\end{example}

\begin{example}
\label{examplenormalbayes}
The general normal model with conjugate normal priors is also quite tractable.  Indeed, for model $X|\theta \sim N_d(\theta, \sigma^2 I_d)$ and prior $\theta \sim N_d(\mu, \tau^2 I_d)$, we have $\theta|x  \sim N_d\big(\mu(x), \tau^2(x) I_d\big)$ with $\hat{\theta}_{\pi}(x) \, = \, \mathbb{E}(\theta|x) \, = \, \frac{\tau^2}{\tau^2 + \sigma^2} x \, + \, \frac{\sigma^2}{\tau^2 + \sigma^2} \mu $ and $\tau^2(x)= \tau^2_0\,=\,\frac{\tau^2 \sigma^2}{\tau^2 + \sigma^2}$.   From this, one obtains

\begin{equation}
\label{lossnormalcasebayes}
L|x \, = \, \|\theta - \mathbb{E}(\theta|x) \|^2 \big| x \sim \tau^2_0
 \, \chi^2_d \,,
\end{equation}
again independent of $x$.   For second-stage squared error loss, one obtains $\hat{L}_{\pi}(x) \, = \, \tau^2_0 \, d$.

 As addressed in Remark \ref{remarklosschoice}, second-stage losses of the form
$\rho(\frac{\hat{L}}{L})$ are desirable alternatives given their scale invariance. Table\ref{table1}, provides Bayesian estimators $\hat{L}_{\pi}$ of $L$ for some choices of $\rho$.  The third column is specific to the normal model context here, while the second column expressions apply more broadly.
\begin{table*}[h]
\caption{Bayesian loss estimators}
\label{table1}
\begin{tabular}{@{}cccccc@{}}
\hline
$\rho(t)$ & $\hat{L}_{\pi}$ (general)
& $\hat{L}_{\pi}$ (conditions)  \\
\hline
$\rho_A(t)  = (t^m-1)^2$    & $\Big(\frac{\mathbb{E}(L^{-m}|x)}{\mathbb{E}(L^{-2m}|x)}\Big)^{1/m}$ & $2 \,\tau^2_0  \Big(\frac{\Gamma(\frac{d}{2}-m)}{\Gamma(\frac{d}{2}-2m)}\Big)^{1/m} \hbox{ for } d > 4m$   \\
$\rho_m(t) = \, t^m - m\log t -1 $    & $ (\mathbb{E}(L^{-m}|x))^{-1/m} $  & $ 2 \tau^2_0 \big(\frac{\Gamma(\frac{d}{2})}{\Gamma(\frac{d}{2} -m)}  \big)^{1/m}  \hbox{ for } d > 2m $  \\
$\rho_B(t) = \,  t + \frac{1}{t} -2$    & $\sqrt{\frac{\mathbb{E}(L|x)}{\mathbb{E}(\frac{1}{L}|x)}} $ & $ \tau^2_0 \sqrt{d(d-2)} \, \hbox{ for } d \geq 3 $   \\
$\rho_C(t) = \,  \big( \log t \big)^2$    & $ e^{\mathbb{E}(\log L|x)} $  & $  2 \tau^2_0 e^{\psi(\frac{d}{2})} \, \hbox{ for } d \geq 1 $  \\
\hline
\end{tabular}
\end{table*}

Various degrees of shrinkage or expansion in comparison to second-stage squared error loss; for which $\hat{L}_{\pi}(X)= d \tau^2_0$; are observable.   For instance with loss $\rho_A$ and $d \geq 5$, we have $\hat{L}_{\pi}(X)= \tau^2_0 (d+2)$ versus $\hat{L}_{\pi}(X)= \tau^2_0 (d-4)$ according to the selections $m=-1$ or $m=1$, respectively.   Shrinkage occurs for $\rho_B$ and $\rho_C$ (see Remark \ref{remarkshrinkage}), while $\hat{L}_{\pi}$ is decreasing as a function of $m$ for both $\rho_A$ and $\rho_m$, with shrinkage iff $m >-1$ for $\rho_m$, and iff $ m > m_0(d)$ for $\rho_A$ with $m_0(d) \in (-1,0)$ such that $\hat{L}_{\pi} = d \tau^2_0$ at $m=m_0(d)$ (see Appendix).

\end{example}
The next examples involve weighted squared error loss as the first-stage which is a typical choice when the model variance varies with $\theta$, such as Poisson and negative binomial.  To this end, consider first-stage loss as weighted squared error loss $L(\theta, \hat{\gamma}) \, = \, \omega(\theta) \, (\hat{\gamma}-\gamma(\theta))^2$ for $X \sim f(x|\theta)$, $\gamma(\theta) \in \mathbb{R}$.   Given a posterior
density for $\theta$, the Bayes estimator of $\gamma(\theta)$ is given, whenever it exists,
 by the familiar
\begin{equation}
\label{Bayesweightedloss}
\hat{\gamma}_{\pi}(x) \, = \, \frac{\mathbb{E}(\gamma(\theta) \, \omega(\theta)|x)}{\mathbb{E}(\omega(\theta)|x)}\,,
\end{equation}
with incurred loss given by $L= \omega(\theta) \big( \hat{\gamma}_{\pi}(x) - \gamma(\theta) \big)^2$.   For second-stage squared error loss $(\hat{L}-L)^2$, we obtain the Bayes estimator
\begin{equation}
\label{Lhatweighted}
\hat{L}_{\pi}(x) \, = \, \mathbb{E}(L|x) \, = \, \mathbb{E}(\gamma(\theta)^2 \omega(\theta)|x) \, - \, \frac{\mathbb{E}^2(\gamma(\theta) \omega(\theta)|x)}{\mathbb{E}(\omega(\theta)|x)}\,,
\end{equation}
as long as $\mathbb{E}(L^2|x)$ exists.

\subsubsection{Poisson distribution}
\label{subsectionPoisson}
Consider a Poisson model $X|\theta \sim \hbox{Poisson}(\theta)$ with a Gamma prior $\theta \sim \hbox{G}(a,b)$ (density proportional to $\theta^{a-1} e^{-\theta b} \, \mathbb{I}_{(0,\infty)}(\theta)$ throughout the manuscript),  and the estimation of
$\gamma(\theta)=\theta$ with normalized squared error loss $\frac{(\hat{\theta}-\theta)^2}{\theta}$.
The set-up leads to $\theta|x \sim \hbox{Ga}(a+x, 1+b)$, and Bayes estimator
\begin{equation}
\label{bayespoisson}
\hat{\gamma}_{\pi}(x) \, = \, \frac{1}{\mathbb{E}(1/\theta|x)} \, = \, \frac{a+x-1}{1+b}\,,
\end{equation}
for $a>1, b>0$.  For the case $(a,b)=(1,0)$, i.e., the uniform prior density on
$(0,\infty)$,   the (generalized) Bayes estimator is also given by (\ref{bayespoisson}), i.e., $\hat{\gamma}_{\pi}(X)=X$, and, moreover, is the unique minimax estimator of $\theta$.

The incurred loss by the Bayes estimator (\ref{bayespoisson}) becomes
\begin{equation}
\label{lossPoisson}
L= \frac{(\frac{a+x-1}{1+b}-\theta)^2}{\theta}\,,
\end{equation}
and the Bayes estimator $\hat{L}_{\pi}$ in (\ref{Lhatweighted}) becomes
\begin{equation}
\label{lpipoisson}
\hat{L}_{\pi}(x) \, = \, \mathbb{E}(\theta|x) \, - \,  \frac{1}{\mathbb{E}(\frac{1}{\theta}|x)} \, = \, \frac{a+x}{1+b} - \frac{a+x-1}{1+b} \, = \, \frac{1}{1+b}\,,
\end{equation}
provided $a > 2$ as the existence of $\mathbb{E}(L^2|x)$ requires finite $\mathbb{E}(\theta^{-2}|x)$ which in turn necessitates $a>2$.    Observe that estimate (\ref{lpipoisson}) is independent of $x$ and of the hyperparameter $a$.

\begin{remark}
\label{remarkLelePoisson}
\begin{enumerate}
\item[ (I)]
Under the above set-up, \cite{Lele1993} established the admissibility of the posterior expectation $\hat{L}_0(X) \, = \, \mathbb{E}(L|X)=1$ as an estimator of $L= \frac{(X-\theta)^2}{\theta}$ under squared error loss $(\hat{L}-L)^2$.   Another interesting property of
$\hat{L}_0(X)$ is that of unbiasedness, as can be seen by the risk calculation $\mathbb{E}(L|\theta)=1$.
\item[ (II)]The frequentist risk of $\hat{L}_{0}(X)$ is given by
\begin{equation}
\nonumber    R(\theta, \hat{L}_{0}) \, = \, \mathbb{E} (\hat{L}_{0}(X)-L)^2 \, = \, \mathbb{V}(L|\theta) \, = \, \frac{\mathbb{E}(X-\theta)^4}{\theta^2} \, - \, 1\, = \,2 + \, \frac{1}{\theta}\,,
\end{equation}
using the fact the fourth central moment of a Poisson distribution with mean $\theta$ is given by $\mathbb{E}(X-\theta)^4 \, = \, \theta (1+3\theta)$.
Observe that the supremum risk is equal to $ + \infty$ which is not conducive to the property of minimaxity.  However, $\hat{L}_0(X)$ is (unique) minimax for estimating $L$ under weighted second-stage loss $\frac{\theta}{2\theta+1} \, (\hat{L}-L)^2$ because it remains admissible, it has constant risk, and such estimators are necessarily minimax.

\item[ (III)]   For $a>2$ and $b>0$, the estimators $\hat{L}_{\pi}(X)$ of $L$   given in (\ref{lpipoisson}) are proper Bayes and admissible since the corresponding integrated Bayes risks are finite.  The finiteness can be justified by the fact that   $\int_0^{\infty}  R(\theta, \hat{L}_{\pi_0}) \, \pi(\theta) \, d\theta \, \leq \, \int_0^{\infty}  R(\theta, \hat{L}_0) \, \pi(\theta) \, d\theta = \, \int_0^{\infty} (2 + \, \frac{1}{\theta}) \, \pi(\theta) d\theta \, = \, 2+ \frac{b}{a-1} < \infty$.

  \end{enumerate}
\end{remark}
\subsubsection{Poisson distribution (multivariate case)}
\label{subsectionPoissonmulti}
As a multivariate extension, consider $X=(X_1, \ldots, X_d)$ with $X_i \sim \hbox{Poisson}(\theta_i)$ independent, the first-stage loss $L=L(\theta, \hat{\theta}) \, = \, \sum_{i=1}^d \frac{(\hat{\theta_i} - \theta_i)^2}{\theta_i}$
for estimating $\theta=(\theta_1, \ldots, \theta_d)$ based on $\hat{\theta}=(\hat{\theta}_1, \ldots, \hat{\theta}_d)$, and second-stage squared error loss.
Proceeding as for the case $d=1$, we have for a given prior $\pi$ the Bayes estimators (whenever well defined):
\begin{equation}
\label{estimatorthetapoissonmulti}
\hat{\theta}_{\pi,i}(x) \,  = \, \big(\mathbb{E} ( \theta_i^{-1}|x)\big)^{-1}\,, i=1, \ldots, d, \\
\end{equation}
\begin{equation}
\label{estimatorLpoissonmulti}
\hbox{and }  \hat{L}_{\pi}(x) \,  =  \, \mathbb{E}(L|x) \, = \,
\sum_{i=1}^d \mathbb{E}(\theta_i|x) \, - \, \sum_{i=1}^d \big(\mathbb{E}(\theta_i^{-1}|x)\big)^{-1}\,.
\end{equation}
As an example, a familiar prior specification choice for $\pi$ (e.g.,  \cite{Clevenson1975}) brings into play $S= \sum_{i=1}^d \theta_i$ and $U_i = \frac{\theta_i}{S}$, $i=1, \ldots, d$, and density  $(S,U) \sim h(s) \, \mathbb{I}_{\{1\}}(\sum_i u_i) $, where $h(\cdot)$ is a density on $\mathbb{R}_+$.  With such a choice and setting $Z=\sum_{i=1}^d x_i$ hereafter, one obtains the posterior density representation $U|s,x \sim \hbox{Dirichlet}(x_1+1, \ldots, x_d+1)$ and $h(s|x) \, \propto \,  s^{Z} \, e^{-s} \, h(s)$.
With $U$ and $S$ independently distributed under the posterior and $\hbox{Beta}(x_i+1, Z-x_i + d-1)$ marginals for the $U_i$'s,
the evaluation of (\ref{estimatorthetapoissonmulti}) and (\ref{estimatorLpoissonmulti}) is facilitated and yields a Clevenson-Zidek type estimator of $\theta$ and accompanying loss estimator
\begin{equation}
\nonumber
\hat{\theta}_{\pi}(X) \, = \, \frac{X}{Z + d -1} \, \big(\mathbb{E}
(S^{-1}|X)\big)^{-1}\, , \hbox{ and } \hat{L}_{\pi}(X) \, = \,  \mathbb{E}(S|X) \,
- \, \frac{Z}{Z+ d -1} \, \big(\mathbb{E}(S^{-1}|X)
\big)^{-1}\, .
\end{equation}
For a gamma prior $S \sim \hbox{G}(a,b)$ with $a \geq 1$, we have $S|x \sim \hbox{G}(a+ Z, b+1)$ and the above reduces to
\begin{equation}
\nonumber
\hat{\theta}_{\pi,a,b}(X) \, = \, \frac{X}{Z + d -1} \, \frac{a + Z -1}{b+1}\,, \hbox{ and } \hat{L}_{\pi,a,b}(X) \, = \,   \frac{ d Z \, + \, a (d-1)}{(b+1) (Z + d -1)} \, \, .
\end{equation}
Notice that the univariate $\hat{L}_{\pi}$ given previously in (\ref{lpipoisson}) is recovered from the above for $d=1$, while the case $a=d, b=0$ yields the unbiased estimators $\hat{\theta}_{\pi,d,0}(X)=X$ and $\hat{L}_{\pi, d, 0}(X)=d$.    We point out that Lele (1992, 1993) established: {\bf (i)} in the latter case, the admissibility of $\hat{L}_{\pi}(X)=d$ as an estimator of loss $L(\theta, X)$ for $d=1,2$, and inadmissibility for $d \geq 3$; and {\bf (ii)} the admissibility of $\hat{L}_{\pi,1,0}$ as an estimator of $L(\theta, \hat{\theta}_{\pi,1,0}(x))$ for all $d \geq 1$.   As in Remark \ref{remarkLelePoisson} for the bivariate case, we point out that $\hat{L}_{\pi, d, 0}(X)=2$ has frequentist risk equal to $4 + \frac{1}{\theta_1} + \frac{1}{\theta_2}$, infinite supremum risk, and that it is minimax for weighted squared error second-stage loss  $\frac{(\hat{L}-L)^2}{4 + \frac{1}{\theta_1} + \frac{1}{\theta_2}}$.

\begin{remark}
\label{remarkPoissonmultivariate}
In the specific situation where $S \sim  \hbox{G}(d,b)$, one verifies that the above prior reduces to independently distributed $\theta_i \sim \hbox{G}(1,b)$ for $i=1, \ldots, d$.  Since the $X_i$'s are also independently distributed given the $\theta_i$'s, the multivariate Bayesian estimation problem reduces to the juxtaposition of $d$ independent univariate problems as analyzed in the previous section.    For instance, expressions (\ref{bayespoisson}) and (\ref{lpipoisson}) applied to the components $\theta_i$ lead to the above $\hat{\theta}_{\pi,d,b}(X)$ and $\hat{L}_{\pi, d, b}(X)$, and the same remains true for the improper proper choice with $b=0$.
  \end{remark}

\subsubsection{Negative binomial distribution}
\label{sectionnb}
Consider a negative binomial model $X \sim \hbox{NB}(r,\theta)$ such that
\begin{equation}
\label{densitynb}
\mathbb{P}(X=x|\theta) \, = \, \frac{(r)_x}{x!}  \, (\frac{r}{\theta+r})^r \, (\frac{\theta}{\theta+r})^x \, \mathbb{I}_{\mathbb{N}}(x)\,,
\end{equation}
with $r>0$, $\mathbb{E}(X|\theta)\,=\, \theta > 0$.   We study pairs
$(\hat{\theta}_{\pi}, \hat{L}_{\pi})$ for a class of Beta type II priors for $\theta$ which are conjugate and defined more generally as follows.

\begin{definition}
\label{definitionBetaII}
A Beta type II distribution, denoted as $Y \sim B2(a,b,\sigma)$ with $a,b,\sigma>0$ has density of the form
\begin{equation}
\nonumber
f(y)\, = \,   \sigma^b   \, \frac{\Gamma(a+b)}{\Gamma(a) \, \Gamma(b)} \, \frac{y^{a-1}}{(\sigma+y)^{a+b}} \, \mathbb{I}_{(0,\infty)}(y)\,.
\end{equation}
The following identity, which is readily verified, will be particularly useful.
\end{definition}

\begin{lemma}\label{lemmabetaIIno}  For $Y \sim B2(a,b,\sigma)$, $\gamma_1 > -a$, and $\gamma_2 > \gamma_1 -b$, we have
\begin{equation}
\label{lemmabetaII}
\mathbb{E} \big(\frac{Y^{\gamma_1}}{(\sigma+Y)^{\gamma_2)}} \big) \, = \, \frac{(a)_{\gamma_1}}{(a+b)_{\gamma_2}}  \,  \frac{(b)_{\gamma_2-\gamma_1}}{\sigma^{\gamma_2 - \gamma_1}} \,,
\end{equation}
where, for $\alpha>0$ and $\alpha+m>0$, $(\alpha)_m$  is the Pochhmamer symbol representing the quantity  $(\alpha)_m \, = \, \frac{\Gamma(\alpha+m)}{\Gamma(\alpha)}$.
\end{lemma}
It is simple to verify the following (e.g., \cite{Ferguson1968}, page 96).

\begin{lemma}
For $X|\theta  \sim \hbox{NB}(r,\theta)$ with prior $\theta \sim B2(a,b,r)$, the posterior distribution is $\theta|x \sim  B2(a+x,b+r,r)$.
\end{lemma}
Now with such a prior, for estimating $\theta$, since $\mathbb{V}(X|\theta) \, = \, \theta (\theta+r)/r$, under normalized squared error loss
\begin{equation}
\label{lossnb}
L(\theta, \hat{\theta}) \, = \,  \frac{(\hat{\theta}- \theta)^2}{\theta (\theta+r)}\,,
\end{equation}
the Bayes estimate of $\theta$ may be derived from (\ref{Bayesweightedloss}) and (\ref{lemmabetaII}) as
\begin{equation}
\nonumber
\hat{\theta}_{\pi}(x) \,=\, \frac{\mathbb{E}\big( \frac{1}{(\theta+r)}|x  \big)}{\mathbb{E}\big( \frac{1}{\theta (\theta+r)}|x  \big)} \, = \,  r \, \frac{a+x-1}{b+r+1}\,,
\end{equation}
for $a+x > 1$.  For $a=1$ and $x=0$, a direct evaluation yields $\hat{\theta}_{\pi}(0)=0$, which matches the above extended to $x=0$.
The associated loss $L(\theta, \hat{\theta}_{\pi}(x))$ has posterior expectation for $a>1$ equal to
\begin{eqnarray*}
\hat{L}_{\pi}(x) \, & = &\,  \mathbb{E}\big( \frac{\theta}{\theta+r}|x  \big) \, - \,
\frac{\mathbb{E}^2\big( \frac{1}{\theta+r}|x  \big)}{\mathbb{E}\big( \frac{1}{\theta (\theta+r)}|x  \big)} \, \\
\, & = & \,  \frac{a+x}{a+b+x+r} \, - \, \big(\frac{b+r}{b+r+1} \big) \,  \big(\frac{a+x-1}{a+b+x+r} \big) \, = \, \frac{1}{b+r+1}\,,
\end{eqnarray*}
making use of (\ref{Lhatweighted}) and (\ref{lemmabetaII}).   Interestingly, the estimator does not depend on $a$ and is constant as a function of $x$ and this property will play a key role for the minimax findings of Section \ref{sectionrukhin}.    We conclude by pointing out that the above applies to the improper prior $\theta \sim B2(1,0,r)$ yielding the generalized Bayes estimator $\hat{\theta}_0(x)\,=\,  rx/(r+1)$.  It is known (e.g., Ferguson, 1968) that the estimator $\hat{\theta}_0$ is minimax with minimax risk equal to $1/(r+1)$.

\subsection{Posterior distributions for loss that do not depend on $x$}
\label{secionposteriorfreeofx}
There are a good number of instances; some of which have appeared in the literature; where both the posterior distribution of loss $L(\theta, \hat{\gamma})$ and (consequently) the Bayes estimate with respect to loss $W(L, \hat{L})$ are free of the observed $x$.   Such a property is particularly interesting and will play a critical role for the minimax implications in Section \ref{sectionminimax}.   We describe situations where such a property arises and collect some examples.   The situations correspond to similar scenarios mathematically and relate specifically to cases where the posterior density admits: (I) a location invariant, (II) a scale invariant, or (III) a location-scale invariant structure.

\begin{theorem}
\label{theoremlocation}
Suppose that the posterior density for $\gamma(\theta)=\theta$ is, for all $x$, location invariant of the form  $\theta|x \sim f(\theta - \mu(x))$ and that the first-stage loss for estimating $\theta$ is location invariant, i.e., of the form $L(\theta, \hat{\theta})= \beta(\hat{\theta}-\theta)$; $\beta: \mathbb{R}^d \to \mathbb{R}_+$.   Then,
\begin{enumerate}
\item[ {\bf (a)} ]   the Bayes estimator $\hat{\theta}_{\pi}(X)$, whenever it exists, is of the form $\hat{\theta}_{\pi}(x) \, = \, \mu(x) \, + \, k$, $k$ being a constant;
\item[ {\bf (b)} ]  the posterior distribution of the loss  $\beta(\hat{\theta}_{\pi}(x)-\theta)$ is free of $x$;
\item[ {\bf (c)} ]  the Bayes estimator $\hat{L}_{\pi}(x)$ of the loss $\beta(\hat{\theta}_{\pi}(x)-\theta)$ with respect to second-stage loss   $W(L, \hat{L})$ is, whenever it exists, free of $x$ .
\end{enumerate}
\end{theorem}
\begin{proof}
Part {\bf (c)} follows from part {\bf (b)}.    Now, observe that
\begin{eqnarray*}
\inf_{\hat{\theta} \in \mathbb{R}^d} \, \mathbb{E} \big\{ \beta(\hat{\theta} - \theta)|x   \big\} & = & \, \inf_{\hat{\theta} \in \mathbb{R}^d} \, \mathbb{E} \big\{ \beta\big(\hat{\theta} - \mu(x) - (\theta - \mu(x)\big)|x \\
\, & = & \, \inf_{\hat{\alpha} \in \mathbb{R}^d} \, \mathbb{E} \big\{ \beta(\hat{\alpha} - Z)|x \big\} \\
\, & = & \, \, \mathbb{E} \big\{ \beta(\hat{\alpha}_{\pi}(x) - Z)|x \big\} \,,
\end{eqnarray*}
for all $x$, with $\hat{\alpha}=\hat{\theta} - \mu(x)$, $\hat{\alpha}_{\pi}(x)\,=\,\hat{\theta}_{\pi}(x) - \mu(x)$, and $Z|x =^d \theta - \mu(x)|x$.  Part  {\bf (a)} follows since the distribution of $Z|x$ does not depend on $x$ and hence the minimizing $\hat{\alpha}$ (i.e., $\hat{\alpha}_{\pi}(x))$ does not depend on $x$.   Finally, part {\bf (b)} follows since
$$  \beta \big(\hat{\theta}_{\pi}(x) - \theta  \big)|x  \, =^d \, \beta\big( k -Z  \big)|x \, $$
is free of $x$.
\end{proof}

We pursue with similar findings as described in (II) and (III) above
\begin{theorem}
\label{theoremscale}
Suppose that the posterior density for $\theta$ is, for all $x$, scale invariant of the form  $\theta|x \sim \frac{1}{\sigma(x)} f(\frac{\theta}{\sigma(x)})$, $\theta \in \mathbb{R}_+$, and that the first-stage loss for estimating $\theta$ is scale invariant, i.e., of the form $L(\theta, \hat{\theta})= \rho(\frac{\hat{\theta}}{\theta})$; $\rho: \mathbb{R} \to \mathbb{R}_+$.  Then,
\begin{enumerate}
\item[ {\bf (a)} ]   the Bayes estimator $\hat{\theta}_{\pi}(X)$, whenever it exists, is of the form $\hat{\theta}_{\pi}(x) \, = \, k \sigma(x) \, $, $k$ being a constant;
\item[ {\bf (b)} ]  the posterior distribution of the loss  $\rho(\frac{\hat{\theta}_{\pi}(x)}{\theta})$ is free of $x$;
\item[ {\bf (c)} ]  the Bayes estimator $\hat{L}_{\pi}(x)$ of the loss   $\rho(\frac{\hat{\theta}_{\pi}(x)}{\theta})$ with respect to second-stage loss   $W(L, \hat{L})$ is, whenever it exists, free of $x$ .
\end{enumerate}
\end{theorem}
\begin{proof}
 A similar development to the proof of Theorem \ref{theoremlocation} establishes the results.
\end{proof}

The next result inspired initially by the context of estimation of a multivariate normal mean with unknown covariance matrix (see Example \ref{examplenormalgamma}) is presented in a more general setting.

\begin{theorem}
\label{theoremlocationscale}
Suppose that the posterior density for $\theta=(\theta_1, \theta_2)$ is, for all $x$, location-scale invariant of the form
\begin{equation}
\label{densitylocationscale}
\theta|x \, \sim \, \frac{1}{\sigma(x) \, \theta_2^d} \, \, f \big(\frac{\theta_1 - \mu(x)}{\theta_2}, \frac{\theta_2}{\sigma(x)}  \big)\,,
\end{equation}
with $\theta_1 \in \mathbb{R}^d, \theta_2 \in \mathbb{R}_+$,
and that the first-stage loss for estimating $\theta_1$ is location-scale invariant, i.e., of the form $L(\theta, \hat{\theta}_1)= \rho\big(\frac{\hat{\theta}_1-\theta_1}{\theta_2}\big)$.   Then,
\begin{enumerate}
\item[ {\bf (a)} ]   the Bayes estimator $\hat{\theta}_{1,\pi}(X)$, whenever it exists, is of the form $\hat{\theta}_{1, \pi}(x) \, = \, \mu(x) \, + \, k \, \sigma(x)$, $k$ being a constant;
\item[ {\bf (b)} ]  the posterior distribution of the loss $L\big(\theta, \hat{\theta}_{1, \pi}(x)\big)$ is free of $x$;
\item[ {\bf (c)} ]  the Bayes estimator $\hat{L}_{\pi}(x)$ of the loss $L\big(\theta, \hat{\theta}_{1, \pi}(x)\big)$  with respect to second-stage loss $W(L, \hat{L})$ is, whenever it exists, free of $x$ .
\end{enumerate}
\end{theorem}
\begin{proof}
Part {\bf (c)} follows from part {\bf (b)}.  Observe that
\begin{equation}
\nonumber   \inf_{\hat{\theta}_1 \in \mathbb{R}^d} \mathbb{E} \big\{ \rho(\frac{\hat{\theta}_1 - \theta_1}{\theta_2})\big| x   \big\} \, = \,
\inf_{\alpha \in \mathbb{R}^d} \mathbb{E} \big\{ \rho(\frac{\alpha - Z}{V})\big| x \,\big\}\,,\
\end{equation}
with $\alpha = \frac{\hat{\theta}_1 - \mu(x)}{\sigma(x)}$, $Z|x =^d \frac{\theta_1 - \mu(x)}{\sigma(x)}| x$, and $V|x =^d \frac{\theta_2}{\sigma(x)}|x$.   Since the pair $(Z,V)|x$ has joint density $\frac{1}{v^d} \, f(\frac{z}{v}, v)$ which is free of $x$, the minimizing $\alpha$ is free of $x$ which yields part {\bf (a)}.  Finally for part
{\bf (b)}, observe that $\rho(\frac{\hat{\theta}_{1, \pi}(x) - \theta_1}{\theta_2})\big| x
\, =^d \,  \rho(\frac{k-Z}{V})|x$, which is indeed free of $x$ given the above.
\end{proof}

\subsubsection{Examples}
First examples that come to mind are given by the non-informative prior density choices:
\begin{enumerate}
\item[ {\bf (i)}]  $\pi(\theta)=1$ for the location model density $X|\theta \sim f_0(x-\theta)$, $x, \theta \in \mathbb{R}^d$, with $f(t)=f_0(-t)$ and $\mu(x)=x$ (Theorem \ref{theoremlocation});

\item[ {\bf (ii)}]  $\pi(\theta)=\frac{1}{\theta}$ for the scale model density
$X|\theta  \sim \frac{1}{\theta} \, f_1(\frac{x}{\theta})$, $x, \theta \ \in \mathbb{R}_+$, with $f(u)=\frac{1}{u^2} f_1(\frac{1}{u})$ and $\sigma(x)=x$ (Theorem \ref{theoremscale});

\item[ {\bf (iii)}]  $\pi(\theta) \, = \frac{1}{\theta_2} \, \mathbb{I}_{(0,\infty)}(\theta_2) \, \mathbb{I}_{\mathbb{R}^d}(\theta_1) $ for $X=(X_1, X_2)|\theta  \sim \frac{1}{\theta_2^{d+1}} \, f_{0,1}\big(\frac{x_1-\theta_1}{\theta_2}, \frac{x_2}{\theta_2}   \big)$ with $f(u,v) \, = \, \frac{1}{v^2} f_{0,1}(-u, \frac{1}{v})$, $\mu(x)\, = x_1\,$, $\sigma(x)=x_2$ (Theorem \ref{theoremlocationscale})\,.
\end{enumerate}
Applications of the above theorems are however not limited to such improper priors and we pursue with further proper prior examples.

\begin{example}  (Multivariate normal model with known covariance)
\label{examplenormal}
Theorem \ref{theoremlocation} applies for the normal model set-up of
Example \ref{examplenormalbayes} since the posterior distribution is of the form
$f(\theta-\mu(x))$.    For instance, under second-stage squared error loss, identity (\ref{lossnormalcasebayes}) and $\hat{L}_{\pi}(X)\, = \, \tau^2_0 d$ are illustrative of parts {\bf (b)} and {\bf (c)} of the theorem.

Theorem \ref{theoremlocation} applies as well to many other first-stage and second-stage losses, such as $L^q$ and reflected normal first-stage losses $\beta(t)= \|t\|^q$ and $1-e^{-c \, \|t\|^2}$, with $c>0$; and second-stage losses of the form $\rho(\frac{\hat{L}}{L})$ such as those of Remark \ref{remarklosschoice}.  Finally, as previously mentioned, the above observations apply for the improper prior density $\pi(\theta)=1$ with corresponding expressions obtained by taking $\tau^2=+\infty$.
\end{example}

\begin{example} (A Gamma model)
\label{examplegamma}
Gamma distributed sufficient statistics appear in many contexts and we consider here $X|\theta \sim  \hbox{G}(\alpha, \theta)$ with a Gamma distributed prior $\theta \sim  \hbox{G}(a, b)$, which results in the scale invariant form of Theorem \ref{theoremscale} with $f$ a $\hbox{G}(\alpha+a, 1)$ density and $\sigma(x) \, = \, (x+b)^{-1}$.  Therefore, Theorem \ref{theoremscale} applies for scale invariant losses $\rho(\frac{\hat{\theta}}{\theta})$ as those referred to in Remark \ref{remarklosschoice}. As an illustration, take entropy-type losses of the form with $\rho_m(t)\,=\, t^m  \, - \, m \log (t) \, - 1$, $m \neq 0$,  for which one obtains for $m < a+  \alpha$ the Bayes estimator
\begin{equation}
\nonumber  \hat{\theta}_{\pi}(x) \, = \, \big\{\mathbb{E}(\theta^{-m}|x  \big\}^{-1/m}
\,= \, k \, \sigma(x) \,,
\end{equation}
 with $k= \big\{ \frac{\Gamma(a+\alpha)}{\Gamma(a+\alpha-m)}   \big\}^{1/m}$, and the posterior distribution of the loss  $\rho_m(\frac{\hat{\theta}_{\pi}(x)}{\theta})$  free of $x$ and matching that of $\rho_m(Z^{-1})$ (or equivalently $\rho_{-m}(Z)$ ) with $Z \sim \hbox{G}(a+\alpha, k)$.   Finally, the Bayes estimate $\hat{L}_{\pi}(x)$ will also be free of $x$ for any second-stage loss.  For the case of squared error loss $W(L,\hat{L})\, = \, (\hat{L}-L)^2$, one obtains
\begin{eqnarray}
\nonumber \hat{L}_{\pi}(x) \, & = & \, \mathbb{E} \big\{\rho_m(\frac{\hat{\theta}_{\pi}(x)}{\theta})|x \big\} \\
\nonumber \, & = & \,  \mathbb{E} \big\{\rho_m(Z^{-1}) \big\} \\
\nonumber \, & = & \,  \mathbb{E} \big(Z^{-m} \, + \, m \mathbb{E} \log Z \, -  \, 1 \big) \\
\label{expressionhatLgamma} \, & = & \,  m \, \Psi(a+\alpha)  \, + \, \log \big( \frac{\Gamma(a+\alpha-m)}{\Gamma(a+\alpha)}  \big)\,,
\end{eqnarray}
for $m < a+ \alpha$, where $\Psi$ is the Digamma function given by $\Psi(t) \, = \, \frac{d}{d t} \log \Gamma (t)$.  To conclude, as previously mentioned, we point out that the above expressions are applicable for the improper density $\pi_0(\theta) \, = \, \frac{1}{\theta}$ on $(0,\infty)$ by setting $a=b=0$.   Related minimax properties are investigated in Section \ref{subsectiongamma}

\end{example}

\begin{example}  (An exponential location model) \\
\label{exanpleexponential}
Consider $X_1, \ldots, X_n$ i.i.d. from an exponential distribution with location parameter $\theta$ and density $e^{-(t-\theta)} \, \mathbb{I}_{(\theta, \infty)}(t)$ (fixing the scale without loss of generality) with a Gamma prior $\theta \sim \hbox{G}(a,b)$, $a>0$ and $b=n$.   This yields a posterior density of the form $\theta|x \sim \frac{1}{\sigma(x)} \, f(\frac{\theta}{\sigma(x)})$ with $\sigma(x) \, = x_{(1)}=\, \min\{x_1, \ldots, x_n\}$, and $f(u) \, = \, a u^{a-1} \, \mathbb{I}_{(0,1)}(u)$, i.e., a Beta$(a,1)$ density.  Theorem \ref{theoremscale} thus applies for any first-stage scale invariant and second-stage losses.

As an illustration, consider the entropy-type loss $L(\theta, \hat{\theta})\,=\, \frac{\theta}{\hat{\theta}} \, - \, \log(\frac{\theta}{\hat{\theta}}) \, - \, 1$, yielding $\hat{\theta}_{\pi}(x) \, = \, \frac{a}{a+1} \, x_{(1)}$  and loss $L=L\big(\theta,\hat{\theta}_{\pi}(x)\big) \,$ distributed under the posterior as $ \frac{a+1}{a} U \, - \, \log(U) \, - \, \log(1+\frac{1}{a}) \, - \, 1$, which is indeed free of $x$.   Finally, for squared error second-stage loss, the Bayes estimator of $L$ is given by
$$  \hat{L}_{\pi}(x) \, = \, \mathbb{E}(L|x) \, = \,  \frac{1}{a} \, - \, \log(1+\frac{1}{a})\,, $$
since $\mathbb{E}(U) \, = \, \frac{a}{a+1}$ and $\mathbb{E}(\log U) \, = \, - \frac{1}{a}$.
\end{example}

\begin{example} (Multivariate normal model with unknown covariance)
\label{examplenormalgamma}

Based on $X=(X_1, \ldots, X_n)^{\top}$ with $X_i \sim N_d(\mu, \sigma^2 I_d)$ independently distributed components , setting $\theta=(\theta_1, \theta_2)=(\mu,\sigma)$,
consider estimating $\theta_1$ under location scale invariant loss  $L(\theta, \hat{\theta}_1) \, = \, \rho\big(\frac{\hat{\theta}_1-\theta_1}{\theta_2}\big)$, such as the typical case $\rho(t)\,=\, \|t\|^2$.  For this set-up, a sufficient statistic is given by $(\bar{X},S)$ with $\bar{X}\,=\, \frac{1}{n} \, \sum_{i=1}^n X_i$ and
$S = \sum_{i=1}^n \|X_i \, - \, \bar{X}\|^2$.  Furthermore, $\bar{X}$ and $S$ are independently distributed as $\bar{X}|\theta \sim N_d(\theta_1, (\sigma^2/n) I_d)$ and $S|\theta \sim \hbox{G}(k/2, 1/2\sigma^2)$ with $k=(n-1)d$.

Now consider a normal-gamma conjugate prior distribution for $\theta$ such that
\begin{equation}
\nonumber   \theta_1|\theta_2 \sim N_d(\xi, \frac{\theta_2^2}{c} I_d)\, \hbox{ with }  \frac{1}{\theta_2^2} \sim \hbox{G}(a,b)\,,
\end{equation}
denoted $\theta \sim \hbox{NG}(\xi, c, a, b)$, with hyperparameters $\xi \in \mathbb{R}^d, a,b,c>0$.  Calculations lead to the posterior density
\begin{equation}
\nonumber
\theta|x \sim NG\big(\xi(x), c+n, a(x), b(x) \big)\,,
\end{equation}
with $\xi(x) \, = \, \frac{n \bar{x} + c \xi}{n+c}, a(x)\,=\, a + \frac{d+k}{2}, \hbox{ and }   b(x) \, = \, \frac{s+2b+\frac{nc}{n+c} \|\bar{x} - \xi\|^2}{2}\,.$   The corresponding posterior density of $\theta$ can be seen to match form (\ref{densitylocationscale}) with $\xi(x)$ as given, $\sigma(x)\, = \, \sqrt{b(x)}$, and $f$ the joint density of $(V_1, V_2)$ where $V_1$ and $V_2$ are independently distributed as $V_1 \sim N_d(0, \frac{1}{c+n} I_d)$ and $V_2^{-2} \sim \hbox{G}(a+ \frac{d+k}{2}, 1)$.  The last representation is obtained by the transformation $(\theta_1, \theta_2) \to (V_1=\frac{\theta_1-\xi(x)}{\theta_2}, V_2=\frac{\theta_2}{\sigma(x)})$ under the posterior distribution.

Theorem \ref{theoremlocationscale} thus applies for any first-stage location-scale invariant and second-stage losses.  For instance, the familiar weighted squared error loss $L(\theta, \hat{\theta}_1)\, = \, \frac{\|\hat{\theta}_1 - \theta_1 \|^2}{\theta_2^2}$ leads to $\hat{\theta}_{1,\pi}{x}\,=\, \xi(x)$,
\footnote{  This can be seen as follows.
\begin{equation}
\label{bayesestimatorofmu}
\hat{\theta}_{1, \pi}(X) \, = \, \frac{\mathbb{E}\Big(\frac{\theta_1}{\theta_2^2} \, \big| x \, \big)}{\mathbb{E}\big(\frac{1}{\theta_2^2}
\, \big| x \, \big)}\, = \,  \frac{\mathbb{E}^{\theta_2|x} \big(\mathbb{E}\big(\frac{\theta_1}{\theta_2^2} \, \big| x, \theta_2 \, \big)\Big)}{\mathbb{E}\big(\frac{1}{\theta_2^2}
\, \big| x \, \big)}\, \, =   \,   \frac{\mathbb{E}(\frac{\xi(x)}{\theta_2^2}|x)}{\mathbb{E}(\frac{1}{\theta_2^2}|x)} \, = \, \xi(x)\,.
\end{equation}}
and loss $L=\frac{\|\theta_1 - \xi(x) \|^2}{\theta_2^2}$ whose posterior distribution (i.e., $\|V_1\|^2$ with $V_1$ as above) is that of a $\frac{1}{n+c} \, \chi^2_d(0)$ distribution.

\end{example}

\begin{remark}
\label{remarkbernstein-vonMises}
Further potential applications of Theorems \ref{theoremlocation} and \ref{theoremlocationscale} may arise when the posterior distribution can be well approximated by a multivariate normal distribution.   Such a situation occurs with the Bernstein-von Mises theorem and the convergence (under regularity conditions; e.g., \cite{DasGupta2008} for an exposition) of
$\sqrt{n} \big(\theta - \hat{\theta}_{mle}  \big) \big| x $ to a $N_d(0, I_{\theta_0})$, $I_{\theta_0}$ being the Fisher information matrix at the true parameter $\theta_0$, and $\theta_{mle}$ the maximum likelihood estimator.
\end{remark}

\section{Minimax findings for a given loss}
\label{sectionminimax}

In this section, we present different scenarios with loss estimators that are minimax and  therefore a benchmark against which we can assess other loss estimators.  The results are subdivided into two parts:  (i) multivariate normal models with independent components and common variance, with or without a known variance; and (ii)  univariate gamma models.   Throughout, the first-stage estimator and associated loss are given, and the task consists in estimating the loss.  For normal models, a quite general class of first-stage losses which are functions of squared error is considered with various second-stage losses of the form $\rho(\frac{\hat{L}}{L})$, while the analysis for the gamma distribution involves first-stage entropy loss and second-stage squared error loss.  The theoretical results are accompanied by observations and illustrations.

\label{subsectionminimax}

\subsection{Normal models}
\label{subsectionnormalmodels}
We first study models $X \sim N_d(\theta, \sigma^2 I_d)$ with known $\sigma^2$ before addressing the unknown $\sigma^2$ case with an i.i.d. sample.   We consider first-stage losses of the form $\beta\big(\frac{\|\hat{\theta} - \theta \|^2}{\sigma^2}  \big)$
for estimating $\theta$, with $\beta(\cdot)$ absolutely continuous and
strictly increasing on $[0,\infty)$, and the loss
\begin{equation}
\label{lossnormalcase}
L \, = \, \beta\big(\frac{\|X - \theta \|^2}{\sigma^2}  \big)
\end{equation}
incurred by estimator $\hat{\theta}_0(X)=X$.   In this set-up, the first-stage frequentist risk is given by $R(\theta, \hat{\theta}_0) \,=\, \mathbb{E}_{\theta}\big\{ \beta\big(\frac{\|X - \theta \|^2}{\sigma^2} \big)  \big\} \, = \, \mathbb{E} \{\beta(Z)\}$ with $Z \sim \chi^2_d$, and we assume that it is finite.  In the identity case $\beta(t)=t$, the estimator $\hat{\theta}_0(X)$ is best equivariant, minimax and (generalized) Bayes with respect to the uniform prior density $\pi_0(\theta)=1$ (e.g., \cite{FSW2018}).   These properties also hold more general $\beta$, and even for $X \sim f(\|x-\theta\|^2)$ with decreasing $f$ (e.g., \cite{Kubokawa2015}).

As a second-stage loss, consider the entropy-type loss
\begin{equation}
\label{entropyloss}   W(L,\hat{L}) \, = \, \rho_m(\frac{\hat{L}}{L})\,, \hbox{ with }
\rho_m(t) \, = \, t^m \, - \, m \log (t) \, - \, 1\,, m \neq 0,
\end{equation}
and the Bayes estimator $\hat{L}_{\pi_0}(X)$ of loss $L$ with respect to prior density $\pi_0$.   Since the posterior distribution $\frac{\|x - \theta \|^2}{\sigma^2} \big| x$ is $\chi^2_d$  (see Example \ref{examplenormalbayes}) independently of $x$, a direct minimization of the expected posterior loss tells us that
\begin{eqnarray}
\nonumber
\inf_{\hat{L}} \big\{\mathbb{E}\big( \rho_m(\frac{\hat{L}}{L})|x \big)  \big\} \, & = & \,
\mathbb{E}\big( \rho_m(\frac{\hat{L}_{\pi_0}}{L})|x \big) \\
\, & = &
\nonumber
 \,   m \mathbb{E} \big (\log \beta(Z) \big ) \, + \, \log \mathbb{E} \big((\beta(Z))^{-m}  \big)\, \\
 \label{normalminimaxrisk}
\, & = & \bar{R} \, ( say) \,,
\end{eqnarray}
as long as these expectations exist, and for
\begin{equation}
\label{expressionlossminimaxnormal}
\hat{L}_{\pi_0}(x)\, = \, \big\{ \mathbb{E} \big((\beta(Z))^{-m} \big\}^{-1/m}\,, \hbox{ with } Z \sim \chi^2_d.
\end{equation}
Also observe that $\hat{L}_{\pi_0}$ has constant risk $R(\theta, \hat{L}_{\pi_0}) =\bar{R}$.   This also can be seen directly by the equivalence of the frequentist and posterior distributions of $\frac{\|X-\theta\|^2}{\sigma^2}$ (see Example \ref{examplelocationfamily}).

\begin{theorem}
\label{theoremminimaxnormal}
Under the above set-up and assumptions, for estimating the first-stage loss $L$ in  (\ref{lossnormalcase}) under second-stage loss (\ref{entropyloss}), the minimax estimator and risk are given by $\hat{L}_{\pi_0}(X)$ and $\bar{R}$, respectively.
\end{theorem}
\begin{proof}
Since $\hat{L}_{\pi_0}(X)$ has constant risk $\bar{R}=R(\theta, \hat{L}_{\pi_0})$, it suffices to show that $\hat{L}_{\pi_0}(X)$ is an extended Bayes estimator with respect to the sequence of priors $\pi_n \sim N_d(0, n \sigma^2 I_d), n \geq 1$, which is to show that
\begin{equation}
\label{condition}
\lim_{n \to \infty} r_n  \, = \,  \bar{R} \,,
\end{equation}
with $r_n$ the integrated Bayes risk associated with $\pi_n$.
For prior $\pi_n$, we have $\theta|x \sim N_d(\frac{n}{n+1} x, \frac{n \sigma^2}{n+1} I_d)$ implying that $\frac{\|X - \theta \|^2}{\sigma^2}|x \sim \frac{n}{n+1} \, \chi^2_d(\frac{y}{n})$, with $y = \frac{\|x\|^2}{(n+1) \sigma^2}$.  Now, setting $Z_n$ such that $Z_n|y \sim
\frac{n}{n+1} \, \chi^2_d(\frac{y}{n})$ and $y$ as above, the Bayesian optimisation problem for $\pi_n$ results in
\begin{equation}
\nonumber
\hat{L}_{\pi_n}(x)\, = \, \Big\{ \mathbb{E} \Big((\beta(\frac{\|x-\theta\|^2}{\sigma^2}\big)^{-m}|x\Big) \Big\}^{-1/m}\, = \, \Big\{ \mathbb{E} \Big((\beta(Z_n)^{-m}|y\Big) \Big\}^{-1/m}\,,
\end{equation}
and
\begin{eqnarray}
\nonumber
\inf_{\hat{L}} \big\{\mathbb{E}_n\big( \rho_m(\frac{\hat{L}}{L})|x \big)  \big\} \, & = & \,
\mathbb{E}_n\big( \rho_m(\frac{\hat{L}_{\pi_n}}{L})|x \big) \\
\, & = &
\nonumber \,   m \mathbb{E} \big (\log \beta(Z_n) |y \big ) \, + \, \log \mathbb{E} \big((\beta(Z_n))^{-m}|y  \big)\, \\
\, & = &  u_n(y)  \hbox{ say}\,.
\end{eqnarray}
From the above, we have
\begin{equation}
\label{rn}
r_n \, = \, \mathbb{E}\big(u_n(Y)  \big) \,  \hbox{ with Y } = \frac{\|X\|^2}{(n+1) \sigma^2}.
\end{equation}

Now observe that the marginal distribution of $X$ is $N_d(0, (n+1) \, \sigma^2 I_d)$ implying that the marginal distribution of $Y$ is $\chi^2_d$ free of $n$.    Finally, by the dominated convergence theorem with $|u_n(y)| \leq v_n(y)$,
$v_n(y)\,=\, \ \mathbb{E}_n\big( \rho_m(\frac{\hat{L}_{\pi_0}}{L})|x \big)$ for all $n \geq 1$ and $y>0$, and
$\mathbb{E}\big(v_n(Y) \big) \, = \ \bar{R}$ (independently of $n$), and
 an application of Lemma \ref{lemmaAppendix2}, which is stated and proven in the Appendix, we infer that
\begin{eqnarray*}
\lim_{n \to \infty}  r_n  \, & = & \,
\mathbb{E}^Y \Big\{ \lim_{n \to \infty} m \, \mathbb{E} \big (\log \beta(Z_n) |Y \big ) \, + \, \lim_{n \to \infty} \log \mathbb{E} \big((\beta(Z_n))^{-m}|Y  \big)  \Big\} \\
\, & = & \,
\mathbb{E}^Y \Big\{ m \mathbb{E} \big (\log \beta(Z) \big ) \, + \, \log \mathbb{E} \big((\beta(Z))^{-m}  \big)  \Big\} \\
\, & = & \,  m \mathbb{E} \big (\log \beta(Z) \big ) \, + \, \log \mathbb{E} \big((\beta(Z))^{-m}  \big)\,,
\end{eqnarray*}
which is (\ref{condition}) and completes the proof.
\end{proof}

Theorem \ref{theoremminimaxnormal} applies quite generally for various choices of $\beta$ and loss $\rho_m$. and the approach is unified.    A particular interesting case is given by first-stage $L^q$ losses with $\beta(t)=t^q \,, q >0$.    Calculations are easily carried out with the moments of $Z \sim \chi^2_d$ yielding the minimax loss estimator $\hat{L}_{\pi_0}(X) \, = \, 2^q \, \big(\frac{\Gamma(\frac{d}{2})}{\Gamma(\frac{d}{2} - mq)}\big)^{1/m}$ of $L= \frac{\|x-\theta\|^{2q}}{\sigma^{2q}}$ for $m < d/2q$.

  An analogous approach establishes the minimaxity of the generalized Bayes loss estimator $\hat{L}_{\pi_0}$ for various other interesting loss functions of the form $\rho(\frac{\hat{L}}{L})$.  We summarize such findings as follows.

\begin{theorem}
\label{theoremotherrho's}
Consider the set-up of Theorem \ref{theoremminimaxnormal} with its corresponding assumptions, and the problem of estimating the first-stage loss (\ref{lossnormalcase}) under second-stage loss $\rho_j(\frac{\hat{L}}{L}), j=A,B,C$ where $\rho_A(t)=(t^m - 1)^2, \, m \neq 0$, $\rho_B(t)\, = \, \frac{1}{2} (t + \frac{1}{t}-2)$, and $\rho_C(t)\,=\, (\log t)^2$,
then assuming existence and finite risk,  the generalized Bayes estimator $\hat{L}_{\pi_0,j}(X)\,, j=A,B,C,$ with respect to the uniform prior density is minimax, its frequentist risk is constant and matches the minimax risk.
\end{theorem}
\begin{proof}
In each of the three cases, a proof is quite analogous to that of Theorem \ref{theoremminimaxnormal} with estimators $\hat{L}_{\pi_0,j}(X)$, the integrated Bayes risk $r_{n,j}, n \geq 1$ and the minimax risk varying for $j = A, B, C$ as presented in the following Table with $Z_n|Y \sim \frac{n}{n+1} \, \chi^2_d(\frac{Y}{n})$ and $Z \sim \chi^2_d$.

\begin{table*}
\caption{}
\label{table2}
\begin{tabular}{@{}cccccc@{}}
\hline
$\rho_j(\frac{\hat{L}}{L})$ & $\hat{L}_{\pi_0,j}(X)$
& $r_{n,j}\,=\, \mathbb{E} \big( g_j(Y)\big), Y \sim \chi^2_d$  &  Minimax \,risk \\
\hline
$\Big((\frac{\hat{L}}{L})^m -1  \Big)^2$    & $\Big\{\frac{\mathbb{E}\big(\beta^{-m}(Z)\big)}{\mathbb{E}\big(\beta^{-2m}(Z)\big)}\Big\}^{1/m}$ & $1 \, - \, \frac{\mathbb{E}^2 \big(\beta^{-m}(Z_n)|Y\big)}{\mathbb{E}\big(\beta^{-2m}(Z_n)|Y\big)}$ &  $ 1 \, - \, \frac{\mathbb{E}^2 \big(\beta^{-m}(Z)\big)}{\mathbb{E}\big(\beta^{-2m}(Z)\big)}$ \\
$\frac{1}{2} \big(\frac{\hat{L}}{L} + \frac{L}{\hat{L}} -2  \big)$ &  $\; \sqrt{\frac{\mathbb{E}\big(\beta(Z)\big)}{\mathbb{E}\big(\frac{1}{\beta(Z)}\big)}}$ & \; $\sqrt{\mathbb{E}\big(\beta(Z_n)|Y \big) \mathbb{E}\big(\frac{1}{\beta(Z_n)}|Y\big)} \, - \, 1 \,$ &  $\sqrt{\mathbb{E}\big(\beta(Z)\big) \mathbb{E}\big(\frac{1}{\beta(Z)}\big)} \, - \, 1$ \\
$\big(\log (\frac{\hat{L}}{L})\big)^2$ & \; $e^{\mathbb{E}\big(\log \beta(Z)\big)}$  &  \; $\mathbb{V}\big(\log \beta(Z_n)\big|Y\big)$  &  $\mathbb{V}\big(\log \beta(Z)\big)$   \\
 \\
\hline
\end{tabular}
\end{table*}
\end{proof}
\begin{remark}
\label{remarkshrinkage}
\begin{enumerate}
\item[ (I)]
For squared error second-stage loss $(\hat{L}-L)^2$, the generalized Bayes estimator $\hat{L}_{\pi_0}$ of $L \, = \, \beta\big(\frac{\|X - \theta \|^2}{\sigma^2}  \big)$ with respect to the uniform prior density $\pi_0$ is given by $\mathbb{E}(L|x)\,=\, \mathbb{E}\big(\beta(Z)  \big)$ with $Z \sim \chi^2_d$ (as for entropy loss $\rho_{-1}$), assuming finite $\mathbb{E}\big(\beta^2(Z)  \big)$.
The methodology of Theorems \ref{theoremminimaxnormal} and \ref{theoremotherrho's} can be applied to establish the minimaxity of $\hat{L}_{\pi_0}$.  This represents an extension of the identity case $\beta(t)=t$ proven by \cite{Johnstone1988}.

\item[ (II)]
Johnstone established in the identity case and for the squared error second-stage loss the inadmissibility of $\hat{L}_{\pi_0}(X)\,=\, d$ for $d \geq 5$ by producing estimators $\hat{L}$ that dominate
$\hat{L}_{\pi_0}$.  These estimators $\hat{L}$, and others appearing later in the literature, are ``shrinkers'' exploiting a potential defect of $\hat{L}_{\pi_0}$.   In comparison, it can be shown using various applications of Jensen's inequality and the covariance inequality that the Bayes estimators $\hat{L}_{\pi_0}$ in Theorem \ref{theoremotherrho's} for $\rho_B$, $\rho_C$, and for $\rho_A$ with $m >0$, as well as those of Theorem \ref{theoremminimaxnormal} for $m >-1$  are shrinkers in the sense that
$\hat{L}_{\pi_0}(X) <  \mathbb{E}\big(\beta(Z)  \big)$.    In contrast, they are expanders for $m <-1$ for both $\rho_m$ and $\rho_A$.  Such comparisons apply as well beyond the set-up here, for other models and priors, namely for Example \ref{examplenormalbayes}'s general $\hat{L}_{\pi}$ expressions, and in comparison to the benchmark posterior expectation estimator $\mathbb{E}(L|X)$  (see Appendix).  Stronger properties can undoubtedly be established in specific cases, such as the identity case seen in Example \ref{examplenormalbayes}. Finally, we point out for a given loss of Theorem \ref{theoremminimaxnormal}, or Theorem \ref{theoremotherrho's}, that the benchmark unbiased procedure $\hat{L}_0(X)=\mathbb{E}\big(\beta(Z)  \big)$ is dominated in terms of frequentist risk by $\hat{L}_{\pi_0}(X)$ unless these two estimators coincide (e.g., $\rho_{-1}$).

\item[ (III)]  The considerations above also informs us on a ``conservativeness'' criterion for selecting a loss estimator which stipulates that
\begin{equation}
\label{propertyconservative}
\mathbb{E}_{\theta} \hat{L}(X)  \geq \mathbb{E}_{\theta} \, L(\theta, \hat{\gamma}(X))\, \hbox{ for all } \theta,
\end{equation}
(equality being (\ref{definitionunbiased})), put forth by \cite{Brown1978}, \cite{LuBerger1989}, and others.    In our context, such an ``expander'' property does not follow in general, and rather is inherited (or disinherited) by the choice of the second-stage loss.  The adherence to (\ref{propertyconservative}) is rather involved in general, but it thus can be controlled in this Example by the choice of the second-stage loss.
\end{enumerate}
\end{remark}

\subsubsection{Unknown $\sigma^2$}
\label{sectionnormalunknownsigma^2}
For the unknown $\sigma^2$ case, a familiar argument (e.g., \cite{LC1998}) coupled with properties of $\hat{L}_{\pi_0}(X)$ for the known $\sigma^2$ case, leads to minimax findings via the following lemma.

\begin{lemma}
\label{lemmaLC}
For $X=(X_1, \ldots, X_n)^{\top}$ with independently distributed components $X_i \sim N_d(\mu, \sigma^2 I_d)$, $\theta=(\mu, \sigma^2)$, consider estimating the loss   $L= \beta(\frac{\|\bar{x} - \mu\|^2}{\sigma^2/n})$ incurred by $\hat{\gamma}(X)=\bar{X}$ for estimating $\gamma(\theta)=\mu$, with $\beta(\cdot)$ absolutely continuous and strictly increasing as in Section \ref{subsectionnormalmodels}.  Suppose that $\hat{L}_0(X)$ is under second-stage loss $W(L,\hat{L})$ free of $\sigma^2$, minimax, and with constant risk $\bar{R}\,=\, \mathbb{E}_{\theta} \big\{ W(L, \hat{L}_0(X))\big\}$ for estimating $L$ regardless of $\sigma^2$.    Then, $\hat{L}_0(X)$ remains minimax for unknown $\sigma^2$ with minimax risk $\bar{R}$.
\end{lemma}
\begin{proof}
Suppose, in order to establish a contradiction, that there exists another estimator $\hat{L}(X)$ such that
$$ \sup_{\theta} \mathbb{E}_{\theta} \big\{ W(L, \hat{L}(X))\big\} \, \, <
\sup_{\theta} \mathbb{E}_{\theta} \big\{ W(L, \hat{L}_0(X))\big\}\,= \, \bar{R}\,. $$
Then, for fixed $\sigma^2=\sigma^2_0$, we would have
$$ \sup_{\theta=(\mu, \sigma^2_0)} \mathbb{E}_{\theta} \big\{ W(L, \hat{L}(X))\big\} \, \, <
\sup_{\theta} \mathbb{E}_{\theta} \big\{ W(L, \hat{L}_0(X))\big\}\,< \, \bar{R}\,, $$  which would be not possible given the assumed minimax property of $\hat{L}_{0}(X)$ for $\sigma^2=\sigma^2_0$.
\end{proof}
The above coupled with the results of the previous section leads to the following.

\begin{corollary}
\label{corollarynormalunknownsigma}
In the set-up of Lemma \ref{lemmaLC} with second-stage loss (\ref{lossnormalcase}), the loss estimator $\hat{L}_{\pi_0}$ given in (\ref{expressionlossminimaxnormal}) is minimax for estimating $L= \beta(\frac{\|\bar{x} - \mu\|^2}{\sigma^2/n})$.   Furthermore, minimaxity is also achieved by the generalized Bayes estimators $\hat{L}_{\pi_0,i}(X)$ of Theorem \ref{theoremotherrho's} for the corresponding losses $W_i(L,\hat{L})\,, \, i=1,2,3.$
\end{corollary}
\begin{proof}
The results are deduced immediately as consequences of Lemma \ref{lemmaLC}, Theorem \ref{theoremminimaxnormal}, and Theorem \ref{theoremotherrho's}.
\end{proof}

\subsection{Gamma models}
\label{subsectiongamma}
We revisit here the Gamma model $X|\theta \sim \hbox{G}(\alpha, \theta)$ of Example \ref{examplegamma} with first-stage entropy-type loss $\rho_m(\frac{\hat{\theta}}{\theta})$ for estimating $\theta$ and with $m < \alpha/2$.      We consider the loss associated with $\hat{\theta}_{\pi_0}(X) \, = \, \frac{k}{X}$, $k= \big\{ \frac{\Gamma(\alpha)}{\Gamma(\alpha-m)}   \big\}^{1/m}$, which as an estimator of $\theta$, is  generalized Bayes for the improper density $\pi_0(\theta)\,=\, \frac{1}{\theta} \, \mathbb{I}_{(0,\infty)}(\theta)$, as well as minimax with constant risk.

With the above first-stage estimator given, the task becomes to estimate
\begin{equation}
\label{lossgamma}
L \, = \, \big( \frac{k}{\theta x}  \big)^m  \, - \, m \, \log \big( \frac{k}{\theta x}  \big) \, - \, 1\ ,
\end{equation}
and we investigate second-stage squared error loss $(\hat{L}-L)^2$.
We establish below the minimaxity of the Bayes estimator
\begin{equation}
\label{Lhatpi0}
\hat{L}_{\pi_0}(X) \, = \,
 m \, \Psi(\alpha)  \, + \, \log \big( \frac{\Gamma(\alpha-m)}{\Gamma(\alpha)}  \big)\,,
\end{equation}
given in Example \ref{examplegamma} for $a=b=0$.  We will require the following Gamma distribution properties, which are derivable in a straightforward manner, and related frequentist risk expression for $\hat{L}_{\pi_0}$.

\begin{lemma}
\label{lemmagammaidentities}
Let $W \sim \hbox{G}(\xi, \beta)$ and $h > -\xi/2$.  We have:
\begin{enumerate}
\item[ {\bf (a)}] $ \mathbb{V}(W^h) \, = \beta^{-2h} \Big\{\, \frac{\Gamma(\xi+2h)}{\Gamma(\xi)} \, - \, \  \big( \frac{\Gamma(\xi+h)}{\Gamma(\xi)} \big)^2 \Big\} $,
\item[ {\bf (b)}] $ \mathbb{V}(\log W) \, = \, \Psi'(\xi)$,
\item[ {\bf (c)}]  $\hbox{Cov}(W^h, \log W) \, = \, \frac{1}{\beta^h} \, \frac{\Gamma(h+\xi)}{\Gamma(\xi)} \, \big\{ \Psi(\xi+h) - \Psi(\xi)  \big\}$,
\end{enumerate}
$\Gamma$ and $\Psi$ being the Gamma and Digamma functions.
\end{lemma}

\begin{lemma}
\label{lemmariskfunctiongamma}
The estimator $\hat{L}_{\pi_0}(X)$ of $L$ has constant in $\theta$ frequentist risk given by
\begin{align}
\label{riskgammaminimax}
R(\theta, \hat{L}_{\pi_0})  \, = k^{2m} \Big\{\, \frac{\Gamma(\alpha-2m)}{\Gamma(\alpha)} \, - \, \  & \notag\big( \frac{\Gamma(\alpha-m)}{\Gamma(\alpha)} \big)^2 \Big\}\,   + \, m^2 \Psi'(\alpha) \, + \, \\
 &  + \,2m \, k^m \, \frac{\Gamma(\alpha-m)}{\Gamma(\alpha)} \, \big\{ \Psi(\alpha+h) - \Psi(\alpha)  \big\},
\end{align}
\end{lemma}
\begin{proof}
Let $Z \sim \hbox{G}(\alpha, k)$. Since $\frac{X\theta}{k}|\theta \, =^d  Z$ and the second-stage loss is squared error, the frequentist risk is equal to $$\mathbb{V}(L|\theta) \, = \, \mathbb{V} \big(Z^{-m} + m \log Z   \big)\, = \,  \mathbb{V}(Z^{-m})   + m^2 \mathbb{V}(\log Z) \, + 2m \, \hbox{Cov}(Z^{-m}, \log Z)\,. $$
The result then follows by applying Lemma \ref{lemmagammaidentities} to the above for $\xi=\alpha, h=-m$.
\end{proof}

We now proceed with the main result of this section.
\begin{theorem}
\label{theoremminimaxgamma}
For $X \sim \hbox{G}(\alpha,\theta)$ with known $\alpha$ ($\alpha>2m$) and unknown $\theta \in \mathbb{R}_+$, first-stage loss $\rho_m(\frac{\hat{\theta}}{\theta})$, and second-stage squared error loss, the generalized Bayes estimator $\hat{L}_{\pi_0}(X)$ given in (\ref{Lhatpi0}) of $L$ is minimax with minimax risk given by (\ref{riskgammaminimax}).
\end{theorem}
\begin{proof}
 We show that $\hat{L}_{\pi_0}(X)$ is an extended Bayes of estimator of $L$ with respect to the sequence of prior densities $\pi_n:= \hbox{G}(a_n, b_n),$ with $a_n=b_n=\frac{1}{n}, n \geq 1$.  Since the risk $R(\theta, \hat{L}_{\pi_0})\,=\, \bar{R}$ is constant in $\theta$, establishing (\ref{condition}) with $r_n$ the integrated Bayes risk with respect to $\pi_n$ will suffice to prove the above.

We have for a given $n$ and $\pi_n$,
\begin{equation}
r_n \, = \,  \mathbb{E}^X_n   \Big\{ \mathbb{E}_n \big(L - \hat{L}_{\pi_n}(X)\big)^2\,\big| \,X  \big)    \Big\} \, = \, \mathbb{E}^X_n \big(\mathbb{V}_n(L|X) \big)\,,
\end{equation}
where $\hat{L}_{\pi_n}(X) \, = \, \mathbb{E}_n(L|X)$ is the Bayes estimator of $L$, $\mathbb{V}_n(L|X)$ is the posterior variance of $L$, and the expectation
$\mathbb{E}^X_n$ is taken with respect to the marginal distribution of $X$.

Under $\pi_n$, we have  $\theta|x \sim \hbox{G}\big(\alpha+a_n, (x+b_n)\big)$ so that $\frac{\theta x}{k}|x \sim \hbox{G}(\alpha+a_n, \frac{k( x+b_n)}{x})$.   Setting $Y = Y(X) = \frac{k(X+b_n)}{X}$, we can write
\begin{eqnarray}
\nonumber \mathbb{V}_n(L|X) \, & = & \, \mathbb{V}_n \big((\frac{k}{\theta X})^m \, - \, m \log (\frac{k}{\theta X})  \big| X \big) \\
\nonumber \, & = & \, \mathbb{V}_n \big( W^{-m} \, + \, m \log W  \big| Y \big)\,,
\end{eqnarray}
with $W|Y \sim \hbox{G}(\alpha+a_n, Y)$.  Expanding the above variance as in Lemma \ref{lemmariskfunctiongamma} and again making use of Lemma \ref{lemmagammaidentities}, one obtains
\begin{equation}
\label{expressionV_n}
\mathbb{V}_n(L|X) \, = \,   C_n \{Y(X)\}^{2m}   \, + \, m^2  \Psi'(a_n+\alpha) \, + \,  2m \, \{Y(X)\}^{m} D_n\,,
\end{equation}
with
$$C_n \, = \, \frac{\Gamma(a_n+\alpha-2m)}{\Gamma(a_n+\alpha)} \, - \,
\frac{\Gamma^2(a_n+\alpha-m)}{\Gamma^2(a_n+\alpha)}\,,$$
and
$$D_n\, = \, \frac{\Gamma(a_n+\alpha-m)}{\Gamma(a_n+\alpha)} \, \big\{\Psi(a_n+\alpha-m) - \Psi(a_n+\alpha)  \big\} \,.  $$

Now, since $X \sim B2(\alpha, a_n, b_n)$ under $\pi_n$, i.e., a Beta type II distribution as in Definition \ref{definitionBetaII}, we have by virtue of Lemma \ref{lemmabetaIIno}
\begin{equation}
\nonumber
\mathbb{E}_n (Y(X)^{h}) \, = \,   \mathbb{E} \big( (\frac{X}{k( X+b_n)})^{-h} \big) \, = \,  k^h \frac{(\alpha)_{-h}}{(\alpha+a_n)_{-h}}\,  \hbox{ for } h < \alpha\,,
\end{equation}
so that $\lim_{n \to \infty}  \mathbb{E}_n \big((Y(X))^{h}  \big) \, = \, k^h$ for $h=m$ and $h=2m$.  Finally with the above and (\ref{expressionV_n}), we obtain directly
$$  \lim_{n \to \infty} r_n \, = \,   \lim_{n \to \infty} \mathbb{E}^X_n \big\{ \mathbb{V}_n(L|X) \big\} \,= \, \bar{R}\,,
$$
as given in (\ref{riskgammaminimax}), completing the proof.
\end{proof}
\section{Minimaxity for Rukhin-type losses}
\label{sectionrukhin}

Whereas the minimax findings of the previous sections apply to decisions problems that are sequential in nature, i.e., the decision of interest which is that of estimating a loss $L=L\big(\theta, \hat{\gamma}(x)\big) $ is assessed for optimality
after having observed the data $x$, Rukhin's loss in (\ref{lossrukhin}) or (\ref{lossrukhingeneral}) applies to the problem of estimating $(\gamma(\theta), L)$ simultaneously.   Whereas \cite{Rukhin1987, Rukhin1988} investigated questions of admissibility of pairs $(\hat{\gamma}, \hat{L})$, our investigation here pertains to minimaxity.    An estimator $(\hat{\gamma}_m, \hat{L}_m)$ of $(\gamma(\theta), L)$ is defined to be minimax for loss
$W(\theta, \hat{\gamma}, \hat{L})$ if  $ \sup_{\theta}
\{W(\theta, \hat{\gamma}_m, \hat{L}_m)\} \, \leq \,
\sup_{\theta}
\{W(\theta, \hat{\gamma}, \hat{L})\}$ for all $(\hat{\gamma}, \hat{L})$.
Since we will investigate the behaviour of a sequence of Bayes estimators, we
point out that a Bayes estimator $(\hat{\gamma}_{\pi}, \hat{L}_{\pi})$ of $
(\gamma, L)$ under loss $W(\theta, \hat{\gamma}, \hat{L})$ and prior $\pi$ is, whenever it exists,
given by $\hat{L}= \mathbb{E}(L|x)$  and $\hat{\gamma}_{\pi}$ the Bayes point estimator of $\hat{\gamma}(\theta)$ under first-stage loss $L=L\big(\theta, \hat{\gamma}(x)\big) $ (independently of the choice of $h$).

We capitalize on a combination of properties and findings of the previous sections
to establish a minimax result which we frame as follows.

\begin{theorem}
\label{theoremrukhin}
Consider a given model $X \sim f_{\theta}$ and loss $W=W(\theta, \hat{\gamma}, \hat{L})$ as in (\ref{lossrukhingeneral}) for estimating $(\gamma, L)$ with $\gamma = \gamma(\theta)$ and $L\,=\, L(\theta, \hat{\gamma}) \, $.  Suppose that there exists an estimator $(\hat{\gamma}_0, \hat{L}_0)$ and a sequence of proper densities $\{\pi_n; n \geq 1\}$ such that: {\bf (i)} $\hat{\gamma}_0(X)$ is for first-stage loss $L$ an extended Bayes estimator of $\gamma$ with constant risk in $\theta$;
{\bf (ii)} the Bayes estimator $\hat{L}_{\pi_n}(x)$ is for $n \geq 1$ constant as a function of $x$;  and {\bf (iii)}  the estimator $\hat{L}_{0}(x)$ is the limit of $\hat{L}_{\pi_n}(x)$ as $n  \to \infty$.   Then $(\hat{\gamma}_0, \hat{L}_0)$ is minimax.
\end{theorem}
\begin{proof}
Denote $R_W\big(\theta, (\hat{\gamma}, \hat{L}) \big) $
as the frequentist risk of estimator $(\hat{\gamma}, \hat{L})$ under loss
$W(\theta, \hat{\gamma}, \hat{L})$; $R(\theta, \hat{\gamma})$ as the first-stage risk of $\hat{\gamma}$; $\bar{R}\,$ as the constant first-stage risk of $\hat{\gamma}_0$; $r_n$ and $r^W_n$ as the integrated Bayes risks with respect to $\pi_n$ associated with the first-stage loss $L$ and global loss $W$, respectively.  As well, denote the constant values of $\hat{L}_{\pi_n}(x)$ and $\hat{L}_{0}(x)$ as $c_n$ and $c=\lim_{n \to \infty}{c_n}$.

We have
\begin{eqnarray*}
R_W\big(\theta, (\hat{\gamma}_0, \hat{L}_0) \big) \, & = & \, \
\mathbb{E}_{\theta}  \big( W(\theta, \hat{\gamma}_0, \hat{L}_0)   \big) \\
\, & = & \,  h'(c) \bar{R} \, - \, c h'(c) \, + \, h(c) \\
 \, & = &  \bar{R}_W\,,
\end{eqnarray*}
which is constant as a function of $\theta$. To establish the result, it will suffice to show that
\begin{equation}
\label{minimaxconditionrukhin}
\lim_{n \to \infty} r^W_n \, = \,  \bar{R}_W.
\end{equation}
 As above, it is the case that
\begin{equation}
\nonumber
R_W(\theta, (\hat{\gamma}_{\pi_n}, \hat{L}_{\pi_n}) \big) \, =  \,
h'(c_n) \, R(\theta, \hat{\gamma}_{\pi_n}) \, -  \, c_n h'(c_n) \, + \, h(c_n)\,,
\end{equation}
which implies that
\begin{equation}
\nonumber
r^W_n \, = \, h'(c_n) \, r_n \, -  \, c_n h'(c_n) \, + \, h(c_n)\,.
\end{equation}
Finally, condition (\ref{minimaxconditionrukhin}) is verified with the above expressions since, by assumptions, $\hat{\gamma}_0$ is extended Bayes with $\lim_{n \to \infty} r_n \, = \, \bar{R}$ and $\lim_{n \to \infty} c_n=c$.

\end{proof}

Observe that the result is quite general and the minimaxity holds irrespectively of the choice of $h$ in loss function (\ref{lossrukhingeneral}), as is the case for the determination of a Bayes estimator of $(\gamma(\theta), L)$.
The above theorem paves the way for various applications which build on the results contained in the previous sections and we present as a series of examples.
A critical property is the one where the Bayes estimators $\hat{L}_{\pi_n}(x)$ are free of $x$, situations that were expanded on in Sections \ref{sectionbayesian}.

\begin{example}  (Normal model)
\label{examplerukhinnormal}
Theorem \ref{theoremrukhin} applies for $X \sim N_d(\theta, \sigma^2 I_d)$, $\gamma(\theta)=\theta$, simultaneous loss $W$ as in (\ref{lossrukhingeneral}) with first-stage squared error loss $L=\frac{\|\hat{\theta}-\theta\|^2}{\sigma^2}$, with the estimator $(\hat{\gamma}_0(X)=X, \hat{L}_0(X)=d)$ which is generalized Bayes for $(\theta, L)$ and the uniform prior density $\pi(\theta)=1$.  Indeed, with the prior sequence of densities $\theta \sim^{\pi_n} N_d(0, n \sigma^2)$, the minimaxity follows from Theorem \ref{theoremrukhin} since: {\bf (i)}  $\hat{\gamma}_0(X)$ is extended Bayes relative to $\{\pi_n; n \geq 1 \}$ with constant risk equal to $d$,   {\bf (ii)} $\hat{L}_{\pi_n}(x) \, = \, \frac{n d}{n+1}$ is constant as a function of $x$, and converges to {\bf (iii)} converges to $\hat{L}_0(x)=d$. \\
\end{example}

\begin{example}  (Normal model continued)
\label{examplerukhinnormal2}
The above example can be extended to the choice of $L= \beta \big(\frac{\|\hat{\theta}-\theta\|^2}{\sigma^2}  \big) $ with $\beta$ a continuous and strictly increasing function on $\mathbb{R}_+$.  Let $Z \sim \chi^2_d(0)$.  Then, the estimator $(\hat{\gamma}, \hat{L})$ with $\hat{\gamma}_0(x)=x$ and $\hat{L}_0(x) \, = \, \mathbb{E} \big( \beta(Z) \big)$ (as long as the latter is finite) can be shown to be minimax for estimating $(\theta,L)$ under loss $W$.  It is also generalized Bayes for the uniform prior density.

A justification of condition {\bf (i) } of the theorem is as follows.  A result in \cite{Yadegari2017} (Theorem 2.2, page 24), that applies when the posterior distribution of $\theta$ is normal, tells  us that the first-stage Bayes estimator of $\theta$ under loss $L$ and prior $\pi_n$ is given by the posterior mean $\frac{nx}{n+1}$, independently of $\beta$,
the posterior being $\theta|x \sim N_d\big(\frac{nx}{n+1}, (\frac{n\sigma^2}{n+1}) I_d \big)$ under prior $\pi_n$.  It follows from this that the minimum expected posterior loss, equivalently $\hat{L}_{\pi_n}(x)$, is equal to $\mathbb{E} \big( \frac{n}{n+1} \, \beta(Z) \big)$ since $ \big(\frac{\|\hat{\theta}-\theta\|^2}{\sigma^2}  \big)|x \sim \frac{n}{n+1} \chi^2_d(0)$.  Now, since this is free of $x$, one infers that the integrated Bayes risk $r_n$ is equal to the minimum expected posterior loss and thus converges to $\mathbb{E} \big( \beta(Z) \big)$ which can be seen as an application of Lemma \ref{lemmaAppendix2} (for $y=0$).  Since this matches the constant risk of $\hat{\gamma}_0(X)$, we infer that the latter is also extended Bayes, whence condition {\bf (i)} of Theorem \ref{theoremrukhin}.    From the above, we infer have that {\bf (ii)} $\hat{L}_{\pi_n}(x)\,=\, \mathbb{E} \big( \frac{n}{n+1} \, \beta(Z) \big)$ is free of $x$, and which {\bf (iii)} converges to $\hat{L}_0(x)$, establishing the minimaxity. \\
\end{example}

\begin{example}  (Gamma model)
\label{examplerukhingamma}
For a Gamma model $X \sim \hbox{G}(\alpha, \theta)$ (i.e., Example \ref{examplegamma}, we apply Theorem \ref{theoremrukhin} for estimating $\gamma(\theta)=\theta$ and $L(\theta, \hat{\theta})= \big(\frac{\hat{\theta}}{\theta}\big)^m \, - \, m \, log (\frac{\hat{\theta}}{\theta}) \, - \,1$ simultaneously under loss $W$ in (\ref{lossrukhin}), with $m < \alpha$.
We show that the Bayes estimator $(\hat{\theta}_0, \hat{L}_0)$ of $(\theta, L)$ with respect to the improper prior density $\pi(\theta) \, = \, \frac{1}{\theta} \, \mathbb{I}_{(0,\infty)}(\theta)$, given by $\hat{\theta}_0(X) \, = \,
\big\{ \frac{\Gamma(a+\alpha)}{\Gamma(a+\alpha-m)}   \big\}^{1/m} \, \frac{1}{X}
\, $ and  $\hat{L}_0(X) \, = m \, \Psi(\alpha)  \, + \, \log \big( \frac{\Gamma(\alpha-m)}{\Gamma(\alpha)}  \big) $ is minimax.   Indeed, with the sequence of prior $\hbox{G}(\frac{1}{n}, \frac{1}{n})$ densities $\pi_n$ , the minimaxity follows since: {\bf (i)} $\hat{\theta}_0(X)$ can be shown to be extended Bayes with constant risk given by $\hat{L}_0(X)$, {\bf (ii)} the Bayes estimator $\hat{L}_{\pi_n}(x)$ of $L$ is a constant given by (\ref{expressionhatLgamma}) with $a=1/n$, and which {\bf (iii)} converges to $\hat{L}_0$ as $n \to  \infty$.

Theorem \ref{theoremrukhin} also applies for other first-stage losses, such as the familiar scale invariant squared error loss $L(\theta, \hat{\theta})\, = \, (\frac{\hat{\theta}}{\theta} -1 )^2$ with $\alpha>2$.  In this case, a minimax solution is $\hat{\theta}_0(X) \, = \, \frac{\alpha-2}{X}$ and $\hat{L}_{0}(X)\, = \, \frac{1}{\alpha-1}$, and the conditions of the theorem can be verified with the same prior sequence $\{\pi_n\}$ as above with $\hat{L}_{\pi_n}(X) \, = \, \frac{1}{\alpha-1+n^{-1}}$ computable from (\ref{Lhatweighted}).
\end{example}

\begin{example}  (Poisson model)
\label{examplerukhinpoisson}
Theorem \ref{theoremrukhin} leads to the following application for the Poisson models of Section \ref{subsectionPoisson} and \ref{subsectionPoissonmulti}.
With the set-up of Section \ref{subsectionPoissonmulti}, for estimating $(\theta, L)$ under loss (\ref{lossrukhin}), we infer that $\big(\hat{\theta}_0(X), \hat{L}_{\pi}(X)\big)=(X,d)$ is minimax by considering the sequence of priors $\pi_n$ such that $S \sim \hbox{G}(a_n=d, b_n=\frac{1}{n})$.  Indeed for such a sequence,
we may show, namely by using Remark \ref{remarkPoissonmultivariate}, that: {\bf (i)} $\hat{\theta}_0(X)$ is extended Bayes with constant risk given by $d$, {\bf (ii)} the Bayes estimator $\hat{L}_{\pi_n}(x)$ of $L$ is a constant as a function of $x$ given by $\frac{d}{1+\frac{1}{n}}$, and which {\bf (iii)} converges to $\hat{L}_0$ as $n \to  \infty$.
\end{example}

\begin{example}  (Negative binomial model)
\label{examplerukhnb}
We consider the set-up of Section \ref{sectionnb} with $X \sim \hbox{NB}(r,\theta)$ as in (\ref{densitynb}), and the problem of estimating $(\theta,L)$ for simultaneous loss $W$ as in (\ref{lossrukhin}),  with $L=L(\theta, \hat{\theta})$ the weighted squared error loss given in (\ref{lossnb}).   Theorem \ref{theoremrukhin} applies in establishing the minimaxity of $(\hat{\theta}_0, \hat{L}_0)$ with $\hat{\theta}_0(X)= \frac{r X}{r+1}$ and $\hat{L}_0(X) \, = \, \frac{1}{r+1}$.  Indeed with the sequence of $\hbox{B}2(a_n, b_n, r)$ prior densities $\pi_n$ for $\theta$ with $a_n=1$ and $b_n=\frac{1}{n}$, it is the case that:  {\bf (i)} $\hat{\theta}_0(X)$ is extended Bayes with constant risk $\bar{R}= \frac{1}{r+1}$, {\bf (ii)}  $\hat{L}_{\pi_n}(x) \, = \, \frac{1}{r+1+n^{-1}}$ is free of $x$, and {\bf (iii)} converges to $\hat{L}_0(x)$ as $n \to \infty$.
\end{example}

The results above paired with the particular features of the $(\hat{\theta}_0, \hat{L}_0)$ minimax solutions lead to further minimax estimators with the simple observation that $(\hat{\theta}_1, \hat{L}_0)$ dominates $(\hat{\theta}_0, \hat{L}_0)$ under loss (\ref{lossrukhin}) whenever $\hat{\theta}_1$ dominates $\hat{\theta}_0$ under first-stage loss $L(\theta, \hat{\theta})$, given that $\hat{L}_0(X)$ is a constant.   We thus have the following implications for the multivariate normal and Poisson models, for $ d \geq 3$ and $d \geq 2$ respectively, since there exist (many) dominating estimators $\hat{\theta}_1(X)$ of $\hat{\theta}_0(X)=X$.  The same applies for the multivariate normal model with a loss function which is a concave function of squared error loss and $d \geq 4$ (see for instance, \cite{FSW2018}).

\begin{corollary}

\begin{enumerate}
\item[ {\bf (a)}]   For the normal model context of Example \ref{examplerukhinnormal} with $d \geq 3$, an estimator $(\hat{\theta}, \hat{L})$ is minimax for estimating $(\theta, L)$ whenever $\hat{\theta}(X)$ dominates
$\hat{\theta}_0(X)=X$ under first-stage loss $\frac{\|\hat{\theta} - \theta\|^2}{\sigma^2}$;

\item[ {\bf (b)}]   For the normal model context of Example \ref{examplerukhinnormal2} with $d \geq 4$ and concave $\beta$, an estimator $(\hat{\theta}, \hat{L})$ is minimax for estimating $(\theta, L)$ whenever $\hat{\theta}(X)$ dominates
$\hat{\theta}_0(X)=X$ under first-stage loss $\beta\big(\frac{\|\hat{\theta} - \theta\|^2}{\sigma^2}\big)$;

\item[ {\bf (c)}]   For the Poisson model context of Example \ref{examplerukhinpoisson} with $d \geq 2$, an estimator $(\hat{\theta}, \hat{L})$ is minimax for estimating $(\theta, L)$ whenever $\hat{\theta}(X)$ dominates
$\hat{\theta}_0(X)=X$ under first-stage loss $ \sum_{i=1}^d  \frac{(\hat{\theta}_i - \theta_i)^2}{\theta_i}\,$.

\end{enumerate}
\end{corollary}

\section{Concluding remarks}

This paper brings into play original contributions and analyses for loss estimation that culminate with minimax findings for:  {\bf (i)} estimating a first-stage loss $L=L(\theta, \hat{\gamma})$, and for {\bf (ii)} estimating jointly  $(\gamma(\theta), L)$ under a Rukhin-type loss.   Various models and choices of the first and second-stage losses were analysed.    Our work also clarifies the structure of various Bayesian solutions, properties of which become critical for the minimax analyses.

All in all, the optimality properties obtained here serve as a guide on how one can sensibly report on an incurred loss in both situations {\bf (i)} and {\bf (ii)}.   Notwithstanding existing results, related questions of admissibility questions remain unanswered, namely in the context of Example \ref{exampleJohnstone} for different second-stage losses where it would be interesting to revisit the effect of the dimension $d$ in the $d-$variate normal case.

\section*{Acknowledgement}
\'{E}ric Marchand's research is supported in part by the Natural Sciences and Engineering Research Council of Canada.   Christine Allard is grateful to the ISM (Institut des sciences math\'ematiques) for financial support.


\begin{appendices}
\section{}
\subsection{}
The following result was used in Section \ref{subsectionnormalmodels}.

\begin{lemma}
\label{lemmaAppendix2}
Let $Z_n \sim \frac{n}{n+1} \, \chi^2_d(\frac{y}{n})$ for $n \geq 1$, $d \geq 1$,   and fixed $y \geq 0$, and let $g$ be a positive valued function such that $\mathbb{E}\big(g(Z_1) \big) < \infty$.   Then, we have $\lim_{n \to \infty} \mathbb{E}\big(g(Z_n) \big) \, = \,  \mathbb{E}\big(g(Z) \big)$ with
$Z \sim \chi^2_d(0)$.
\end{lemma}
\begin{proof}
Denote $h_n$ and $f_n$ as the density of $Z_n$ and $\frac{n+1}{n}Z_n \sim \chi_d^2(\frac{y}{n})$ respectively.
We seek to apply the dominated convergence theorem and we make use of the $W \sim \chi_{\nu}^2(\lambda)$ density representation
$$\big(\frac{1}{2}\Big)^{\nu/2} \, \frac{w^{\frac{\nu -2}{2}}}{\Gamma(\frac{\nu}{2})} \,\,_0F_1(-;\frac{\nu}{2}; \frac{\lambda w}{4}) \, e^{-(\lambda + w)/2} \,,$$
for $w > 0$, $\nu>0$, $\lambda \geq 0$. From this, we have for all $n \geq 2$ and $z>0$:
$$\begin{array}{rcl} g(z) h_n(z) & = &  g(z) \, f_n (\frac{(n+1) z}{n}) \, \frac{n+1}{n} \\
 & \leq &  \frac{3}{2} \, g(z) f_n(\frac{n+1}{n}z) \, \\
 & = & \frac{3}{2} \, g(z) \, \big(\frac{1}{2}\big)^{d/2}\frac{(\frac{(n+1)z}{n})^{\frac{d -2}{2}}}{\Gamma(\frac{d}{2})}\,\, _0F_1(-;\frac{d}{2}; \frac{n+1}{4n^2}zy) \; e^{-\frac{(y + (n+1)z)}{2n}} \\
 & \leq & \frac{3}{2} \, K \, e^{\frac{y}{2}} \, g(z) \, f_1(z) \,, \, \end{array}\,$$
where $K = \max\{1,(\frac{3}{2})^{\frac{d-2}{2}}\}$. The result then follows by dominated convergence.

\end{proof}

\subsection{}  As a complement to Example \ref{examplenormalbayes}, here are justifications to the effect that the constant loss estimate $\hat{L}_{\pi}$ decreases as a function of $m$ for losses {\bf (i)} and {\bf (ii)}, and that the defined cut-off point $m_0(d)$ takes values between $-1$ and $0$.    We have
\begin{equation}
\nonumber  \log \big(\frac{\hat{L}_{\pi}(x)}{2 \tau^2_0}  \big) \, = \, f(d,m) \, = \,  \frac{1}{m} \, \big\{   \log \Gamma(\frac{d}{2}) \, - \,  \log \Gamma(\frac{d}{2} - m)\big\}\,,
\end{equation}
for {\bf (ii)}, and $\log \big(\frac{\hat{L}_{\pi}(x)}{2 \tau^2_0}  \big) \, = \, f(d-2m,m)$ for {\bf (i)}.   Now, observe that $f(d,m)$ increases in $d$, and decreases in $m$, the former being a consequence of the strict logconvexity of the gamma function, and the latter since
\begin{equation}
\nonumber
\frac{\partial}{\partial m} f(d,m) \, = \, -\frac{1}{m^2} \Big\{ \log \Gamma(\frac{d}{2}) \, - \,  \log \Gamma(\frac{d}{2} - m) \, - m \, \Psi(\frac{d}{2} -m) \Big\}  < 0\,,
\end{equation}
with the inequality due to the ordering $u'(a) < \frac{u(a+b) - u(a)}{b} < u'(a+b)$ for $a,b \in \mathbb{R}_+$ and strictly convex and differentiable functions $u(\cdot)$ on $\mathbb{R}_+$.   The above tells us directly that $\hat{L}_{\pi}$ decreases in $m$ for loss {\bf (ii)}, but it is also the case for loss {\bf (i)} since $f(d-2m_1, m_1) >  f(d-2m_2,m_1) > f(d-2m_2,m_2)$ for $m_1 < m_2 < \frac{d}{4}$.

For the bounds on $m_0(d)$ which apply to loss {\bf (i)}, it suffices to observe that $\hat{L}_{\pi} = \tau^2_0 (d+2)$ for $m=-1$, calculate the limiting value $\hat{L}_{\pi}= \, 2 \tau^2_0 \, e^{\Psi(\frac{d}{2})}$ as $m \to 0$, and then infer that  $\lim_{m \to 0} \hat{L}_{\pi} \leq \tau^2_0 d$ by virtue of
the Digamma function inequality $\Psi(\alpha) < log(\alpha)$ for $\alpha>0$.  More generally,  one shows that $\lim_{m \to 0} \Big(\frac{\mathbb{E}(L^{-m}|x)}{\mathbb{E}(L^{-2m}|x)}\Big)^{1/m} \, = \, e^{\mathbb{E}(\log L|x)}  $, so that the squared log error loss arises as the limiting loss $\big((\frac{\hat{L}}{L})^m-1\big)^2$ when $m \to 0$.

\subsection{}    Here are elements of justification for the stated properties of Remark \ref{remarkshrinkage}.  We make use of the inequalities $\mathbb{E}\big(h(\beta(Z))  \big) \leq h\big(\mathbb{E}( \beta(Z)) \big)$  for concave $h$ (Jensen), and $\mathbb{E} \big(f(\beta(Z)) \, g(\beta(Z))  \big) \leq \mathbb{E} \big(f(\beta(Z)) \big)\mathbb{E} \big(g(\beta(Z))  \big)$ for increasing $f$ and decreasing $g$ (Covariance).  The implications for losses $\rho_m$ and $\rho_C$ follow with Jensen's inequality using $h(t)= t^{-m} \, \hbox{ or } -t^{-m}$
depending on the value of $m$, and $h(t)=\log(t)$, respectively.   The shrinkage for $\rho_A$ with $m \in (0,1)$ follows with the covariance inequality applied for $f(t)=t^m$ and $g(t)=t^{-2m}$, telling us that $\big(\hat{L}_{\pi_0}(x)\big)^m < \mathbb{E}(L^m|x)$, followed by Jensen's inequality applied to $h(t)=t^m$.  The shrinkage that occurs for $\rho_B$ follows from the use of the covariance inequality
with $f(t)=t$ and $g(t)=t^{-1}$.
There remains loss $\rho_A$, the Bayes estimator $\big\{\frac{\mathbb{E}(L^{-m}|x)}{\mathbb{E}(L^{-2m}|x)}\big\}^{1/m}$ and its properties of shrinkage for $m  > 0$, and expansion for $m <-1$, in comparison to the benchmark estimator $\hat{L}_0(X)=\mathbb{E}(L|X)$.   These properties follow directly from the following inequality, which is also of independent interest.

\begin{lemma}
\label{lemmamoments}
The following inequality holds for a positive random variable $T$:
$$\frac{\mathbb{E}(T^{-m})}{\mathbb{E}(T^{-2m})} \leq \big(\mathbb{E}(T)\big)^m \hbox{ for } m>0 \hbox{ and } m \leq -1\,,$$
assuming existence of the above expectations.
 \end{lemma}
\begin{proof}
For a positive real number $N$, we set $\lfloor{N}\rfloor$ and $\{N\}$ as integer and fractional parts defined here as $\lfloor{N}\rfloor \, = \, \sup\{j \in \mathbb{N}: j < N\}$ and $\{N\} = N - \lfloor{N}\rfloor $. The result has been previously shown for $m \in (0,1)$.  For $m \geq 1$, the result follows by applying the covariance inequality $\lfloor{m}\rfloor + 1 $ times as follows
\begin{eqnarray*}
\mathbb{E}(T^{-m}) \leq \mathbb{E}(T^{-m-1}) \, \mathbb{E}(T) \leq \cdots  &  \leq &
 \mathbb{E}(T^{-m-\lfloor{m}\rfloor}) \, \big(\mathbb{E}(T)   \big)^{\lfloor{m}\rfloor}       \\
\,  &  \leq &  \mathbb{E}(T^{-m-\lfloor{m}\rfloor - \{m\}}) \, \big(\mathbb{E}(T)   \big)^{\lfloor{m}\rfloor + \{m\}}  \\
 \,  &  = &  \mathbb{E}(T^{-2 m} ) \, \big(\mathbb{E}(T)   \big)^{m}\,.
\end{eqnarray*}
Similarly for $m<-1$, the inequality follows as
\begin{eqnarray*}
\mathbb{E}(T^{-2m}) \geq \mathbb{E}(T^{-2m-1}) \, \mathbb{E}(T) \geq \cdots  &  \geq &
 \mathbb{E}(T^{-2m-\lfloor{-m}\rfloor}) \, \big(\mathbb{E}(T)   \big)^{\lfloor{-m}\rfloor}       \\
\,  &  \geq &  \mathbb{E}(T^{-2m-\lfloor{-m}\rfloor - \{-m\}}) \, \big(\mathbb{E}(T)   \big)^{\lfloor{-m}\rfloor + \{-m\}}  \\
 \,  &  = &  \mathbb{E}(T^{- m} ) \, \big(\mathbb{E}(T)   \big)^{-m}\,.
\end{eqnarray*}
\end{proof}
\end{appendices}

\end{document}